\documentclass[12pt]{article}
\usepackage{amssymb}
\usepackage{amsmath}
\usepackage{array}
\usepackage[T2A]{fontenc}
\usepackage[cp1251]{inputenc}
\usepackage[russian]{babel}
\usepackage[pdftex]{graphics}

\multlinegap=0em \DeclareFontFamily{T1}{msb}{}
\DeclareFontShape{T1}{msb}{m}{ol}{<5> <6> <7> <8> <9> gen * msbm
<10> <10.95> <12> <14.4> <17.28> <20.74> <24.88> msbm10}{}
\DeclareSymbolFont{AMSb}{T1}{msb}{m}{ol} \multlinegap=0em

\renewcommand{\S}{\mathhexbox278}
\renewcommand{\le}{\operatorname{\leqslant}}

\begin{document}

\begin{flushleft}
11M06
\end{flushleft}

\begin{center}
{\rmfamily\bfseries\normalsize M.A. Korolev}
\end{center}

\begin{center}
{\rmfamily\bfseries\normalsize On the large values of the Riemann
zeta-function on short segments of the critical line\footnote{The
research is supported by Russian Scientific Fund (grant No.
14-11-00433).}}
\end{center}

\vspace{0.5cm}

\fontsize{11}{12pt}\selectfont

\textbf{Abstract.} In this paper, we obtain a series of new
conditional lower bounds for the modulus and the argument of the
Riemann zeta function on very short segments of the critical line,
based on the Riemann hypothesis.

\vspace{0.2cm}


\vspace{0.2cm}

Bibliography: 36 titles.

\fontsize{12}{15pt}\selectfont

\vspace{0.5cm}

\textbf{Introduction}

\vspace{0.5cm}

In 2001\,-2006, A.A.~Karatsuba
\cite{Karatsuba_2001a}-\cite{Karatsuba_2006} obtained a series of
lower estimates for the maximum of modulus of the Riemann zeta
function $\zeta(s)$ in the circles of small radius lying in the
critical strip $0\leqslant\Re s\leqslant 1$, and on very short
segments of the critical line $\Re s = 0.5$. These results gained
further progress in \cite{Garaev_2002}-\cite{Changa_2005}.

In particular, it was proved in \cite{Karatsuba_2006} that the
function
\[
F(T;H)\,=\,\max_{|t-T|\leqslant H}\bigl|\zeta(0.5+it)\bigr|
\]
satisfies the inequality
\begin{equation}\label{lab_00}
F(T;H)\,\geqslant\,\frac{1}{16}\,\exp\biggl\{-\,\frac{5\ln{T}}{6(\pi/\alpha\,-\,1)(\cosh{(\alpha
H)-1})}\biggr\},
\end{equation}
where $\alpha$ is any fixed number, $1\leqslant\alpha<\pi$,
$2\leqslant \alpha H\leqslant$ $\ln\ln{H}-c_{1}$, and $c_{1}>0$ is
some absolute \hfill constant. \hfill Given \hfill $\varepsilon>0$,
\hfill it \hfill follows \hfill from \hfill (\ref{lab_00}) \hfill
that \hfill for \hfill any \hfill $T\geqslant T_{0}(\varepsilon)$
\hfill and \hfill for
\\ $H\geqslant \pi^{-1}(1+\varepsilon)\ln\ln{T}-c_{1}$, the function $F(T;H)$ is
bounded from below by some constant:
\[
F(T;H)\,>\,c_{2}\,=\,\frac{1}{16}\,\exp{\bigl(-\,1.7\varepsilon^{-1}e^{c_{1}}\bigr)}>0.
\]

In \cite{Karatsuba_2006}, A.A.~Karatsuba posed the problem to prove
the inequality $F(T;H)\geqslant 1$ for the values of $H$ which are
essentially smaller than $\ln\ln{T}$, namely, for  $H\geqslant
\ln\ln\ln{T}$.\footnote{\;If $\ln\ln{T}\ll H\leqslant 0.1T$, then
the following estimate of R.~Balasubramanian
\cite{Balasubramanian_1986} holds true:
\[
F(T;H)\,\gg\,\exp{\biggl(\frac{3}{4}\sqrt{\frac{\ln{H}}{\ln\ln{H}}}\biggr)}.
\]
This bound is supposed to be close to the best possible. Thus, the
estimates of $F(T;H)$ for $0<H\ll \ln\ln{T}$ are most interesting in
this topic.}

In this paper, we give a conditional solution of Karatsuba's
problem, based on the Riemann hypothesis. Moreover, we prove that
for arbitrary large fixed number $A\geqslant 1$ there exist positive
constants $T_{0}$ and $c_{0}$ that depend on $A$ and such that for
any $T\geqslant T_{0}$ and $H = (1/\pi)\ln\ln\ln{T}+c_{0}$ the
inequality $F(T;H)\geqslant A$ holds true (see Theorem 1).

The method used here is applicable to the estimation both of the
maxima of the function
\[
\bigl|\zeta(0.5+it)\bigr|\,=\,\exp{\bigl(\ln{\bigl|\zeta(0.5+it)\bigr|}\bigr)}\,=\,\exp{\bigl(\Re
\ln{\zeta(0.5+it)}\bigr)},
\]
and the extremal values of the function
\[
S(t)\,=\,\frac{1}{\pi}\,\arg{\zeta\bigl(0.5+it\bigr)}\,=\,\frac{1}{\pi}\Im\ln{\zeta(0.5+it)}
\]
(for the definition and basic properties of the function $S(t)$,
which is called the argument of the Riemann zeta\,-function on the
critical line, see the survey \cite{Karatsuba_Korolev_2005}).

The estimates of maximum and minimum of the function $S(t)$ on very
short segments of the variation of $t$ hold the significant interest
together with classical estimates of the values $\displaystyle
\max_{T\leqslant t\leqslant 2T}\bigl(\pm S(t)\bigr)$ belonging to
A.~Selberg \cite{Selberg_1946} and K.\,-M.~Tsang \cite{Tsang_1986},
\cite{Tsang_1993}. Thus, the estimates of the form
\[
\max_{|t-T|\leqslant H}\bigl(\pm S(t)\bigr)\,\geqslant\, f(H),
\]
where
\[
f(H)\,=\,\frac{1}{90\pi}\sqrt{\frac{\ln{H}}{\ln\ln{H}}},\quad
(\ln{T})(\ln\ln{T})^{-\,3/2}<H<T
\]
and
\[
f(H)\,=\,\frac{1}{900}\,\frac{\sqrt{\ln{H}}}{\ln\ln{H}},\quad
\sqrt{\ln\ln{T}}\,\leqslant\,H\,\leqslant\,(\ln{T})(\ln\ln{T})^{-\,3/2}
\]
are obtained in \cite{Korolev_2005} and \cite{Boyarinov_2011a},
\cite{Boyarinov_2011b}, respectively.

In this paper, we prove the existence of positive and negative
values of the function $S(t)$ whose moduli exceed 3, on each segment
of length $H = 0.8\ln\ln\ln{t}+c_{0}$ (see Theorems 2\,-4). For
comparison, we note that it appears in the process of calculation of
first 200 billions zeros of $\zeta(s)$ on the critical line
(S.~Wedeniwski \cite{Wedeniwski_2003}, 2003) that
\begin{align*}
& |S(t)|\,<\,1\quad \text{if}\quad 7<t< 280;\\
& |S(t)|\,<\,2\quad \text{if}\quad 7<t< 6\,820\,050;\\
& |S(t)|\,<\,3\quad \text{if}\quad 7<t< 16\,220\,609\,807.
\end{align*}
The first values of $S(t)$ which exceed 3 in modulus, are located in
the neighborhoods of Gram points $t_{n}$ (see \S 4) with indices $n
= 53\,365\,784\,979$ и $n = 67\,976\,501\,145$ and are equal to
$3.0214$ and $-3.2281$, respectively. At present time, no values of
$t$ such that $|S(t)|\geqslant 4$ are known.

Since the function $S(t)$ is ``responsible'' for the irregularity in
the distribution of zeros of $\zeta(s)$, Theorems 3 and 4 imply some
conditional results related the distribution of Gram's intervals
$G_{n} = (t_{n-1},t_{n}]$ which contain an ``abnormal'' (that is,
unequal to 1) number of ordinates of zeros of $\zeta(s)$ (see
Theorems 5, 6).

The paper ends with a proof of Theorem 7 that concerns the
distribution of nonzero values of integer-valued function
$\Delta_{n}$ introduced by A.~Selberg \cite{Selberg_1946b} in
connection with so-called Gram's law.

In this paper, we use the following notations: $\Lambda(n)$ denotes
the von Mangoldt function, which is equal to $\ln{p}$ for prime $p$
and $n = p^{k}$, $k = 1,2,\ldots$, and equal to zero otherwise;
$\Lambda_{1}(n) = \Lambda(n)/\ln{n}$ ($n\geqslant 2$); $\cosh{z} =
(e^{z}+e^{-z})/2$; $K_{a}(z) = \exp{\bigl(-a\cosh{z}\bigr)}$
($a>0$); $\hat{f}$ denotes the Fourier transform of the function
$f$, that is
\[
\hat{f}(u)\,=\,\int_{-\infty}^{+\infty}f(x)e^{-iux}\,dx;
\]
$\|\alpha\| = \min{\bigl(\{\alpha\},1-\{\alpha\}\bigr)}$ is the
distance between $\alpha$ and the closest integer; $p_{1} = 2$,
$p_{2}=3$, $p_{3}=5, \ldots$ are primes indexed in ascending order;
$\Omega(n)$ is the number of prime factors of $n$ counted with
multiplicity; $\theta, \theta_{1}, \theta_{2}, \ldots$ are complex
numbers, different in different formulae, whose absolute value does
not exceed one. All other notations are explained in the text.

\vspace{0.5cm}

\textbf{\S 1. Auxilliary assertions}

\vspace{0.5cm}

In this section, we give some auxilliary lemmas.

\vspace{0.2cm}

\textsc{Lemma 1.} \emph{For any $m\geqslant 1$, the numbers}
\[
1,\;\frac{1}{2\pi}\ln{2},\;\frac{1}{2\pi}\ln{3},\;\frac{1}{2\pi}\ln{5},\quad\ldots,\quad
\frac{1}{2\pi}\ln{p_{m}}
\]
\emph{are linearly independent over the field of rationals.}

\vspace{0.5cm}

\textsc{Proof.} Let's assume the contrary. Then there exist the
integers $k\geqslant 0, k_{1},$ $k_{2},\ldots, k_{m}$ not equal to
zero simultaneously and such that
\[
k\,+\,\frac{k_{1}}{2\pi}\ln{2}\,+\,\frac{k_{2}}{2\pi}\ln{3}\,+\,\ldots\,+\frac{k_{m}}{2\pi}\ln{p_{m}}\,=\,0,
\]
or, which is the same,
\begin{equation}\label{lab_01}
k\,-\,\frac{1}{2\pi}\ln\frac{a}{b}\,=\,0,
\end{equation}
where $a$ and $b$ are coprime integers not equal to one
simultaneously, whose prime factors do not exceed $p_{m}$.
Exponentiating (\ref{lab_01}), we get
\begin{equation}\label{lab_02}
e^{2\pi k}\,=\,\frac{a}{b}.
\end{equation}

If $k = 0$ then (\ref{lab_02}) contradicts the fundamental theorem
of arithmetics. If $k\geqslant 1$ then $e^{\pi}$ appears to be the
root of polynomial $bz^{2k}-a$. This is impossible in view of
transcendence of $e^{\pi}$ (see for example \cite[\S
2.4]{Feldman_1982}). These contradictions prove the lemma.

\vspace{0.5cm}

\textsc{Lemma 2.} \emph{The estimate
$\bigl|\widehat{K}_{a}(\lambda)\bigr|\leqslant \kappa
e^{-b|\lambda|}$ holds for any real $\lambda$ with
\[
\kappa\,=\,\kappa(a,b)\,=\,2\int_{0}^{+\infty}\exp{\bigl(-a(\cos{b})\cosh{u}\bigr)}\,du,
\]
where $b$ is any number with the condition} $0<b<\pi/2$.

\vspace{0.5cm}

The proof of this statement repeats almost verbatim that of Lemma 4
in \cite{Korolev_2005b}.

\vspace{0.5cm}

\textsc{Lemma 3.}\emph{Suppose that $\lambda$ is real and satisfies
the condition $|\lambda|\geqslant a\sqrt{2}$. Then the following
relation holds:}
\[
\widehat{K}_{a}(\lambda)\,=\,\frac{2\sqrt{2\pi}}{\sqrt[4\;]{\lambda^{2\mathstrut}-a^{2}}}\,
\exp{\biggl(-\,\frac{\pi|\lambda|}{2}\biggr)}\bigl(\cos{g_{a}(\lambda)}\,+\,r_{a}(\lambda)\bigr),
\]
\emph{where}
\[
g_{a}(\lambda)\,=\,\sqrt{\lambda^{2}-a^{2}}\,-\,|\lambda|\ln{\biggl(\frac{|\lambda|}{a}\,+\,
\frac{\sqrt{\lambda^{2\mathstrut}-a^{2}}}{a}\biggr)}\,+\,\frac{\pi}{4},\quad
|r_{a}(\lambda)|\,\leqslant\,c_{a}|\lambda|^{-0.1},
\]
\emph{and}
\begin{equation*}
c_{a}\,=\,
\begin{cases}
9.3, & \emph{if}\quad a\geqslant 1/\sqrt{2},\\
8.2a^{-0.4}, & \emph{if}\quad 0<a<1/\sqrt{2}.
\end{cases}
\end{equation*}

\vspace{0.5cm}

\textsc{Proof.} Without loss of generality, we assume that
$\lambda>0$.  Let's take an arbitrary $R>1$ and denote by
$\Gamma_{R}$ the contour of rectangle with vertices at the points
$\pm R$, $\pm R - \pi i/2$, traversed counterclockwise. The
application of Cauchy's residue theorem yields
\[
\int_{\Gamma_{R}}K_{a}(z)e^{-i\lambda z}\,dz\,=\,\sum\limits_{k
=1}^{4}I_{k}\,=\,0,
\]
where $I_{1}$, $I_{3}$ are integrals along the upper and lower sides
of contour and $I_{2}, I_{4}$ are integrals over lateral sides.

Further, it is easy to note that
 \begin{multline*}
-\, I_{1}\,=\,\int_{-R}^{R}K_{a}(u)e^{-i\lambda u}\,du,\\
 I_{3}\,=\,\int_{-R}^{R}K_{a}\biggl(u\,-\,\frac{\pi i}{2}\biggr)e^{-i\lambda\bigl(u\,-\,
 \frac{\scriptstyle \pi i}{\scriptstyle 2\mathstrut}\bigr)}\,du\,=\,e^{-\,\frac{\scriptstyle \pi \lambda}{\scriptstyle 2\mathstrut}}
 \int_{-R}^{R}e^{i\varphi_{a}(u)}\,du,
 \end{multline*}
where $\varphi_{a}(u)=a\sinh{u}\,-\,\lambda u$. Let us put $z =
R-\pi it/2$, where $0\leqslant t\leqslant 1$. Since
$|K_{a}(z)|=e^{-a\cosh{(R)}\cos{(\pi t/2)}}$, we get:
 \begin{multline*}
 |I_{4}|\,\leqslant\,\frac{\pi}{2}\int_{0}^{1}e^{-a\cosh{(R)}\cos{(\pi t/2)}}\,dt\,=\,
 \frac{\pi}{2}\int_{0}^{1}e^{-a\cosh{(R)}\sin{(\pi t/2)}}\,dt\,\leqslant\,\frac{\pi}{2}\int_{0}^{1}e^{-at\cosh{(R)}}\,dt
 \,\leqslant\\
 \leqslant\,\frac{\pi}{2a\cosh{R}}.
 \end{multline*}
The same bound is valid for the integral $I_{2}$. Hence,
\[
\int_{-R}^{R}K_{a}(u)e^{-i\lambda
u}\,du\,=\,e^{-\,\frac{\scriptstyle \pi \lambda}{\scriptstyle
2\mathstrut}}
\int_{-R}^{R}e^{i\varphi_{a}(u)}\,du\,+\,\frac{\pi\theta}{a\cosh{R}}.
\]
Letting $R$ tend to infinity, we obtain:
\[
\widehat{K}_{a}(\lambda)\,=\,e^{-\,\frac{\scriptstyle \pi
\lambda}{\scriptstyle
2\mathstrut}}\int_{-\infty}^{+\infty}e^{i\varphi_{a}(u)}\,du\,
=\,2e^{-\,\frac{\scriptstyle \pi \lambda}{\scriptstyle
2\mathstrut}}\Re{j_{a}(\lambda)},\quad
j_{a}(\lambda)\,=\,\int_{0}^{+\infty}e^{i\varphi_{a}(u)}\,du.
\]
The derivative $\varphi_{a}'(u)$ has a unique zero on the ray of
integration at a point
\[
u_{a}\,=\,\text{arccosh}\,\frac{\lambda}{a}\,=\,\ln{\biggl(\frac{\lambda}{a}\,+\,\sqrt{\frac{\lambda^{2}}{a^{2}}\,-\,1}\biggr)}.
\]
Setting $u = u_{a}+v$, where $-u_{a}\leqslant v<+\infty$ and noting
that
\[
\varphi_{a}(u)\,=\,a\bigl(\sinh{u_{a}}\cosh{v}\,+\,\cosh{u_{a}}\sinh{v}\bigr)\,-\,\lambda(u_{a}+v)\,=\,-\lambda
u_{a}\,+\,\lambda\psi_{a}(v),
\]
where $\psi_{a}(v) = \alpha\cosh{v}+\sinh{v}-v$, $\alpha =
\sqrt{1-(a/\lambda)^{2\mathstrut}}$, we find that
\[
j_{a}(\lambda)\,=\,e^{-i\lambda
u_{a}}\int_{-u_{a}}^{+\infty}e^{i\lambda \psi_{a}(v)}\,dv.
\]
Suppose that $\delta$ satisfies the condition
$0<\delta<\min{\bigl(1,u_{a},\lambda^{-1/3}\bigr)}$. Then we
represent $j_{a}(\lambda)$ as the sum
\[
e^{-i\lambda
u_{a}}\biggl(\;\int_{-\delta}^{\delta}\,+\,\int_{-u_{a}}^{-\delta}\,+\,\int_{\delta}^{+\infty}\biggr)e^{i\lambda
\psi_{a}(v)}\,dv\,=\, e^{-i\lambda
u_{a}}\bigl(j_{1}+j_{2}+j_{3}\bigr).
\]
We have
\[
\psi_{a}(v)\,=\,\psi_{a}(0)\,+\,\psi_{a}'(0)v\,+\,\psi_{a}''(0)\,\frac{v^{2}}{2}\,+\,\psi_{a}^{(3)}(\xi)\,\frac{v^{3}}{6}
\]
for $|v|\leqslant \delta$, where $\xi$ lies between $0$ and $v$.
Since
\[
\psi_{a}'(v)\,=\,\alpha\sinh{v}+\cosh{v}-1,\quad
\psi_{a}''(v)\,=\,\alpha\cosh{v}+\sinh{v},\quad
\psi_{a}^{(3)}(v)\,=\,\alpha\sinh{v}+\cosh{v},
\]
then $\psi_{a}(0) = \psi_{a}''(0) = \alpha$, $\psi_{a}'(0) = 0$, and
\[
\bigl|\psi_{a}^{(3)}(\xi)\bigr|\,=\,\bigl|\alpha\sinh{\xi}\,+\,\cosh{\xi}\bigr|\,\leqslant\,\sinh|\xi|\,+\,\cosh{\xi}\,=\,e^{|\xi|}\,\leqslant\,e^{\delta}\,<\,e.
\]
Hence,
\[
\lambda\psi_{a}(v)\,=\,\mu\,+\,\mu\,\frac{v^{2}}{2}\,+\,e\lambda\,\frac{\theta
v^{3}}{6},\quad
\mu\,=\,\alpha\lambda\,=\,\sqrt{\lambda^{2\mathstrut}-a^{2}}.
\]
Let us define $\varrho(v)$ by the relation
$\exp{\bigl(ie\theta\lambda v^{3}/6\bigr)}=1+\varrho(v)$. Thus we
get
\begin{multline*}
\bigl|\varrho(v)\bigr|\,=\,\biggl|\frac{ie\lambda}{6}\,\theta
v^{3}\,+\,\frac{1}{2!}\biggl(\frac{ie\lambda}{6}\,\theta
v^{3}\biggr)^{2}\,+\,\frac{1}{3!}\biggl(\frac{ie\lambda}{6}\,\theta
v^{3}\biggr)^{3}+\ldots\biggr|\,\leqslant\\
\leqslant\,\frac{e\lambda}{6}\,|v|^{3}\biggl(1\,+\,\frac{1}{2!}\frac{e}{6}\,+\,\frac{1}{3!}\biggl(\frac{e}{6}\biggr)^{2}+\ldots\biggr)\,=\,
\bigl(e^{e/6}-1\bigr)\lambda|v|^{3}\,<\,\frac{3\lambda}{5}|v|^{3}.
\end{multline*}
Therefore,
\begin{multline*}
j_{1}\,=\,\int_{-\delta}^{\,\delta}\exp{\biggl(i\mu + \frac{i\mu
v^{2}}{2}\biggr)}\bigl(1\,+\,\varrho(v)\bigr)\,dv\,=\,e^{i\mu}\int_{-\delta}^{\,\delta}\exp{\biggl(\frac{i\mu
v^{2}}{2}\biggr)}\,dv\,+\,2\theta_{1}\int_{0}^{\delta}\frac{3\lambda}{5}v^{3}\,dv\,=\\
=\,2e^{i\mu}\int_{0}^{\delta}\exp{\biggl(\frac{i\mu
v^{2}}{2}\biggr)}\,dv\,+\,\frac{3\theta_{1}}{10}\,\lambda\delta^{4}\,=\,e^{i\mu}\sqrt{\frac{2}{\mu}}
\int_{0}^{\frac{\scriptstyle\mu}{\scriptstyle
2\mathstrut}\delta^{2}}\frac{e^{iw}\,dw}{\sqrt{w}}
\,+\,\frac{3\theta_{1}}{10}\,\lambda\delta^{4}.
\end{multline*}
Replacing the last integral by improper one and noting that
\[
\int_{0}^{+\infty}\frac{e^{iw}\,dw}{\sqrt{w}}\,=\,e^{\pi
i/4}\sqrt{\pi},\quad
\biggl|\int_{u}^{+\infty}\frac{e^{iw}\,dw}{\sqrt{w}}\biggr|\,\leqslant\,\frac{2}{\sqrt{u}},
\]
we find that
\begin{multline*}
j_{1}\,=\,e^{i\mu}\sqrt{\frac{2}{\mu}}\biggl(\sqrt{\pi}e^{\pi
i/4}\,+\,\frac{2\theta_{2}\sqrt{2}}{\sqrt{\mu\delta^{2\mathstrut}}}\biggr)\,+\,\frac{3\theta_{1}}{10}\,\lambda\delta^{4}\,=\,\sqrt{\frac{2\pi}{\mu}}\,e^{i(\mu
+ \pi
i/4)}\,+\,\theta_{3}\biggl(\frac{4}{\mu\delta}\,+\,\frac{3\lambda\delta^{4}}{10}\biggr)
\end{multline*}
for any $u>0$. Further, the integration by parts in $j_{2}$ yields:
\[
j_{2}\,=\,\frac{1}{i\lambda}\biggl(\frac{e^{i\lambda\psi_{a}(-\delta)}}{\psi_{a}'(-\delta)}\,-\,
\frac{e^{i\lambda\psi_{a}(-u_{a})}}{\psi_{a}'(-u_{a})}\,-\,\int_{-u_{a}}^{-\delta}e^{i\lambda\psi_{a}(v)}d\,\frac{1}{\psi_{a}'(v)}\biggr)
\]
and hence
\[
|j_{2}|\,\leqslant\,\frac{1}{\lambda}\biggl(\frac{1}{|\psi_{a}'(-\delta)|}\,+\,\frac{1}{|\psi_{a}'(-u_{a})|}\,+\,
\int_{-u_{a}}^{-\delta}\biggl|d\,\frac{1}{\psi_{a}'(v)}\biggr|\biggr).
\]
Since
\[
\alpha\,=\,\frac{\sqrt{\lambda^{2\mathstrut}-a^{2}}}{\lambda}\,=\,\frac{\sinh{u_{a}}}{\cosh{u_{a}}}\,=\,\tanh{u_{a}},
\]
then the derivative $\psi_{a}''(v) =
\cosh{v}\bigl(\alpha+\tanh{v}\bigr)$ is positive for $v>-u_{a}$.
Thus, the function $1/\psi_{a}'(v)$ decreases for $v>-u_{a}$. Hence,
\begin{multline*}
|j_{2}|\,\leqslant\,\frac{1}{\lambda}\biggl(\frac{1}{|\psi_{a}'(-\delta)|}\,+\,\frac{1}{|\psi_{a}'(-u_{a})|}\,-\,
\int_{-u_{a}}^{-\delta}d\,\frac{1}{\psi_{a}'(v)}\biggr)\,=\\
=\,\frac{1}{\lambda}\biggl(\frac{1}{|\psi_{a}'(-\delta)|}\,+\,\frac{1}{|\psi_{a}'(-u_{a})|}\,-\,
\frac{1}{\psi_{a}'(-\delta)}\,+\,\frac{1}{\psi_{a}'(-u_{a})}\biggr).
\end{multline*}
Since $\psi_{a}'(0) = 0$, then $\psi_{a}'(v)<0$ for negative $v$ and
therefore
\[
j_{2}\,\leqslant\,\frac{2}{\lambda|\psi_{a}'(-\delta)|}.
\]
Further, we have
\[
|\psi_{a}'(-\delta)|\,=\,\bigl|\alpha\sinh{\delta}\,-\,\cosh{\delta}\,+\,1\bigr|\,=\,2\sinh{\frac{\delta}{2}}\,\biggl|\alpha\cosh{\frac{\delta}{2}}
\,-\,\sinh{\frac{\delta}{2}}\biggr|\,>\,\delta\biggl|\alpha\cosh{\frac{\delta}{2}}
\,-\,\sinh{\frac{\delta}{2}}\biggr|.
\]
Since $\lambda\geqslant a\sqrt{2}$, then $\alpha\geqslant
1/\sqrt{2}$ and hence
\begin{multline*}
\alpha\cosh{\frac{\delta}{2}}
\,-\,\sinh{\frac{\delta}{2}}\,\geqslant\,\frac{1}{\sqrt{2}}\cosh{\frac{\delta}{2}}
\,-\,\sinh{\frac{\delta}{2}}\,\geqslant\,\frac{1}{\sqrt{2}}\biggl(1\,+\,\frac{1}{2!}\biggl(\frac{\delta}{2}\biggr)^{2}
\,+\,\frac{1}{4!}\biggl(\frac{\delta}{2}\biggr)^{4}\,+\,\ldots\biggr)
,-\\
-\,\biggl(\frac{\delta}{2}\,+\,\frac{1}{3!}\biggl(\frac{\delta}{2}\biggr)^{3}\,+\,\frac{1}{5!}\biggl(\frac{\delta}{2}\biggr)^{5}\,+\,\ldots\biggr)\,>\,
\frac{1}{\sqrt{2}}\,-\,\frac{\delta}{2}\,>\,\frac{1}{5}.
\end{multline*}
Finally we get:
\[
|\psi_{a}'(-\delta)|\,>\,\frac{\delta}{5},\quad
|j_{2}|\,<\,\frac{10}{\lambda\delta}\,<\,\frac{10}{\mu\delta}.
\]
The proof of the inequality $|j_{3}|\leqslant
2\bigl(\lambda\psi_{a}'(\delta)\bigr)^{-1}$ is just the same. By the
relations
$\psi_{a}'(\delta)\,=\,\alpha\sinh{\delta}\,+\,\cosh{\delta}\,-\,1\,>\,\alpha\delta\,\geqslant\,\delta/\sqrt{2}$,
it implies that
\[
|j_{3}|\,\leqslant\,\frac{2\sqrt{2}}{\lambda\delta}\,<\,\frac{3}{\mu\delta}.
\]
Therefore,
\[
j_{1}+j_{2}+j_{3}\,=\,\sqrt{\frac{2\pi}{\mu}}\,e^{i(\mu+\pi
/4)}\,+\,r_{1},
\]
where
\[
|r_{1}|\,\leqslant\,\frac{4}{\mu\delta}\,+\,\frac{3\lambda\delta^{4}}{10}\,
+\,\frac{10}{\mu\delta}\,+\,\frac{3}{\mu\delta}\,=\,\frac{17}{\mu\delta}\,+\,\frac{3\lambda\delta^{4}}{10}.
\]
Thus we conclude that
\[
j_{a}(\lambda)\,=\,\sqrt{\frac{2\pi}{\mu}}\,e^{i(\mu + \pi/4-\lambda
u_{a})}\bigl(1\,+\,r_{2}\bigr),
\]
where
\[
|r_{2}|\,\leqslant\,\sqrt{\frac{\mu}{2\pi}}\biggl(\frac{17}{\mu\delta}\,+\,\frac{3\lambda\delta^{4}}{10}\biggr)\,\leqslant\,
\frac{1}{\sqrt{\pi}}\biggl(\frac{17}{\sqrt[4\;]{2}\,\delta\sqrt{\lambda}}\,+\,\frac{3\lambda^{3/2}\delta^{4}}{10\sqrt{2}}\,\biggr).
\]

If $a\sqrt{2}\geqslant 1$, we put $\delta = (7/8)\lambda^{-2/5}$.
Since $\lambda\geqslant a\sqrt{2}\geqslant 1$, the inequalities
$\delta < 1$, $\delta<\lambda^{-1/3}$ are obvious. Moreover,
\[
u_{a}\,=\,\ln{\biggl(\frac{\lambda}{a}\,+\,\sqrt{\biggl(\frac{\lambda}{a}\biggr)^{2}\,-\,1}\biggr)}\,\geqslant\,\ln{(\sqrt{2}+1)}\,>\,\frac{7}{8}\,\geqslant\,\delta,
\]
and hence $\delta<\min{(1,\lambda^{-1/3}, u_{a})}$. Thus, we have in
this case:
\[
|r_{2}|\,\leqslant\,\frac{1}{\sqrt{\pi}}\biggl(\frac{8\cdot
17}{7\sqrt[4\;]{2}}\,+\,\frac{3}{10\sqrt{2}}\biggl(\frac{7}{8}\biggr)^{\!
4}\biggr)\lambda^{-1/10}\,<\,9.3\lambda^{-0.1}.
\]
If $a\sqrt{2}<1$ then we put $\delta = (a/\lambda)^{2/5}$. Then the
inequality $\lambda\geqslant a\sqrt{2}$ implies that
\[
\delta \leqslant (1/\sqrt{2})^{2/5}<1,\quad
a^{6}\,<\,a\biggl(\frac{1}{\sqrt{2}}\biggr)^{5}\,=\,\frac{a\sqrt{2}}{8}\,\leqslant\,\frac{\lambda}{8}\,<\,\lambda,
\]
and $a^{2/5}<\lambda^{1/15} = \lambda^{2/5-1/3}$. Thus, $\delta <
\lambda^{-1/3}$. Finally, since the inequality
$x^{-2/5}<\ln{\bigl(x+\sqrt{x^{2}-1}\bigr)}$ holds for any
$x\geqslant\sqrt{2}$, we find  $\delta < u_{a}$. Therefore, in this
case, the inequality $\delta<\min{(1,\lambda^{-1/3}, u_{a})}$ is
also valid. Thus we obtain
\[
|r_{2}|\,\leqslant\,\frac{1}{\sqrt{\pi}}\biggl(\frac{17}{\sqrt[4\;]{2}}\,+\,\frac{3a^{2}}{10\sqrt{2}}\biggr)a^{-2/5}\lambda^{-1/10}\,
<\,8.2a^{-0.4}\lambda^{-0.1}.
\]

Finally we get
\begin{multline*}
\widehat{K}_{a}(\lambda)\,=\,2e^{-\pi\lambda/2}\sqrt{\frac{2\pi}{\mu}}\,\Re\biggl(e^{i(\mu-\lambda
u_{a}+\pi/4)}\bigl(1\,+\,r_{2}\bigr)\biggr)\,=\\
=\,2\sqrt{\frac{2\pi}{\mu}}\,e^{-\pi\lambda/2}\bigl(\cos{(\mu-\lambda
u_{a}+\pi/4)}\,+\,r\bigr),
\end{multline*}
where $|r|\leqslant c_{a}\lambda^{-0.1}$ is such that $c_{a} = 9.3$
for $a\sqrt{2}\geqslant 1$ and $c_{a} = 8.2a^{-0.4}$ for
$0<a\sqrt{2}<1$. The lemma is proved.

\vspace{0.2cm}

\textsc{Corollary.} \emph{Under the conditions of Lemma} 3,
\emph{the following inequality holds}
\[
\bigl|\widehat{K}_{a}(\lambda)\bigr|\,<\,\kappa_{a}\,\frac{e^{\,-\,\pi|\lambda|/2}}{\sqrt{|\lambda|}},
\]
\emph{where $\kappa_{a} = 61.5$ for $a\sqrt{2}\geqslant 1$ and
$\kappa_{a} = 54.1a^{-0.4}$ for} $0<a\sqrt{2}<1$.

\vspace{0.2cm}

\textsc{Proof.} The inequality of Lemma 3 together with the
condition $|\lambda|\geqslant a\sqrt{2}$ imply that
\[
\bigl|\widehat{K}_{a}(\lambda)\bigr|\,<\,\frac{2\sqrt{2\pi}}{\sqrt{|\lambda|}}\,
\frac{e^{\,-\,\pi|\lambda|/2}}{\sqrt[4\;]{1-(a/\lambda)^{2\mathstrut}}}\,(1+r)\,\leqslant\,
\frac{2^{7/4}\sqrt{\pi}}{\sqrt{|\lambda|}}\,e^{\,-\,\pi|\lambda|/2}(1+r),
\]
where $r = c_{a}|\lambda|^{-1/10}$. Using the above expressions for
$c_{a}$, we get the desired bound.

\vspace{0.5cm}

\textsc{Lemma 4.} \emph{Suppose that the function $f(z)$ is
analytical in the strip $|\!\Im z|\leqslant 0.5+\alpha$, where it
satisfies the inequality $|f(z)|\leqslant c(|z|+1)^{-(1+\beta)}$
with some positive $\beta$ and $c$. Then the identity}
 \begin{multline}\label{lab_03}
 \int_{-\infty}^{+\infty}f(u)\ln{\zeta\bigl(0.5+i(t+u)\bigr)}\,du\,=\,\sum\limits_{n = 2}^{+\infty}\frac{\Lambda_{1}(n)}{\sqrt{n}}\,n^{-it}\hat{f}(\ln{n})\,+\\
 +\,2\pi\biggl(\;\sum\limits_{\beta>0.5}\int_{0}^{\beta-0.5}f(\gamma-t-iv)\,dv\,-\,\int_{0}^{0.5}f(-t-iv)\,dv\biggr),
 \end{multline}
\emph{holds for any $t$, where $\varrho=\beta+i\gamma$ in the last sum
runs through all complex zeros of $\zeta(s)$ to the right from the
critical line.}

 \vspace{0.5cm}

This assertion goes back to A.~Selberg (see for example \cite[Lemma
16]{Selberg_1946}). In \cite[Ch. II, \S 2]{Karatsuba_Korolev_2005},
\cite{Tsang_1986}, there are some variants of this lemma, where
$f(z)$ satisfies slightly different conditions. These proofs can be
easily adopted to the case under considering.

\vspace{0.5cm}

\textsc{Lemma 5.} \emph{If the Riemann hypothesis is true then the
relation}
\begin{multline}\label{lab_04}
\int_{-\infty}^{+\infty}K_{a}(\pi
u)\ln{\zeta\bigl(0.5+i(t+u)\bigr)}\,du\,=\,\frac{1}{\pi}\sum\limits_{n
=
2}^{+\infty}\frac{\Lambda_{1}(n)}{\sqrt{n}}\,n^{-it}\widehat{K}_{a}\biggl(\frac{\ln{n}}{\pi}\biggr)\,-\\
-\,2\pi\int_{0}^{0.5}K_{a}(\pi t + \pi i v)dv
\end{multline}
\emph{holds for any real $t$.}

\vspace{0.5cm}

\textsc{Proof.} We take an arbitrary $\delta$ such that
$0<\delta<10^{-6}$ and set $z = x+iy$, $f(z) =
K_{a}((\pi-\delta)z)$, $\alpha = \delta/(4\pi)$. Since the
inequalities
\[
\cos{(\pi-\delta)y}\,\geqslant\,\cos{\bigl\{(\pi-\delta)(0.5+\alpha)\bigr\}}\,>\,\sin{\frac{\delta}{4}}\,\geqslant\,2\alpha,
\]
hold for any $y$ such that $|y|\leqslant 0.5+\alpha$, then we have
\[
|f(z)|\,=\,e^{-a\cosh{(\pi-\delta)x}\cos{(\pi-\delta)y}}\,\leqslant\,e^{-2a\alpha\cosh{(\pi-\delta)x}}\,\leqslant\,c(|z|+1)^{-(1+\beta)}.
\]
for a suitable constants $\beta = \beta(\alpha)$, $c = c(\alpha)$
and for any $x$.

The application of Lemma 4 yields:
\begin{multline}\label{lab_05}
\int_{-\infty}^{+\infty}K_{a}\bigl((\pi-\delta)u\bigr)\ln{\zeta\bigl(0.5+i(t+u)\bigr)}\,du\,=\\
=\,\frac{1}{\pi-\delta}\sum\limits_{n =
2}^{+\infty}\frac{\Lambda_{1}(n)}{\sqrt{n}}\,n^{-it}\widehat{K}_{a}\biggl(\frac{\ln{n}}{\pi-\delta}\biggr)\,-
\,2\pi\int_{0}^{0.5}K_{a}\bigl((\pi-\delta)(t+iv)\bigr)\,dv.
\end{multline}
Let us take
\[
N\,=\,\biggl[\frac{1}{\delta^{2}}\biggl(\ln\frac{1}{\delta}\biggr)^{\!-1}\biggr]\,+\,1
\]
and suppose $\delta$ to be so small that $N>N_{0}=e^{\pi
a\sqrt{2}}$. Now we split the sum in (\ref{lab_05}) to the sums
$C_{1}, C_{2}$ and $C_{3}$ corresponding to the intervals $n>N$,
$N_{0}<n\leqslant N$ и $n\leqslant N_{0}$, respectively. Using the
Corollary of Lemma 3 with $\lambda = (1/\pi)\ln{n}\geqslant
a\sqrt{2}$, we obtain
\[
|C_{1}|\,\leqslant\,\frac{1}{\pi-\delta}\sum\limits_{n>N}\frac{\Lambda_{1}(n)}{\sqrt{n}}\,61.5\sqrt{\frac{\pi-\delta}{\ln{n}}}\,
\exp{\biggl(-\,\frac{\pi}{2}\,\frac{\ln{n}}{\pi-\delta}\biggr)}\,\leqslant\,\frac{61.5}{\sqrt{\pi-\delta}}
\sum\limits_{n>N}\frac{\Lambda(n)}{n(\ln{n})^{3/2}}.
\]
The application of Abel's summation formula together with the bound
\begin{equation}\label{lab_06}
\psi(u)\,=\,\sum\limits_{n\leqslant
u}\Lambda(n)\,\leqslant\,c_{1}u,\quad c_{1}\,=\,1.03883
\end{equation}
(see \cite[Th. 12]{Rosser_Schoenfeld_1962}), which is valid for any
$u>0$, yields:
\begin{multline*}
\sum\limits_{n>N}\frac{\Lambda(n)}{n(\ln{n})^{3/2}}\,=\,-\int_{N}^{+\infty}\bigl(\psi(u)\,-\,\psi(N)\bigr)\,
d\,\frac{1}{(\ln{u})^{3/2}}\,\leqslant\,-c_{1}\int_{N}^{+\infty}u\,d\,\frac{1}{(\ln{u})^{3/2}}\,=\\
=\,c_{1}\biggl(\frac{2}{\sqrt{\ln{N}\mathstrut}}\,+\,\frac{1}{(\ln{N})^{3/2}}\biggr).
\end{multline*}
Using the inequalities $\ln{N}\geqslant\ln{\bigl(1/\delta\bigr)}$
and $0<\delta<10^{-6}$, we get the estimate
\[
|C_{1}|\,\leqslant\,\frac{123c_{1}}{\sqrt{\pi-\delta}}\,\frac{1}{\sqrt{\ln{(1/\delta)}}}\biggl(1\,+\,\frac{1}{2\ln{(1/\delta)}}\biggr)
\,<\,\frac{75}{\sqrt{\ln{(1/\delta)}}}.
\]
Similarly,
\[
\biggl|\frac{1}{\pi}\sum\limits_{n>N}\frac{\Lambda_{1}(n)}{\sqrt{n}}\,n^{-it}\widehat{K}_{a}\biggl(\frac{\ln{n}}{\pi}\biggr)\biggr|\,<\,\frac{74.9}{\sqrt{\ln{(1/\delta)}}}.
\]
Thus we get:
\[
C_{1}\,=\,\frac{1}{\pi}\sum\limits_{n>N}\frac{\Lambda_{1}(n)}{\sqrt{n}}\,n^{-it}\widehat{K}_{a}\biggl(\frac{\ln{n}}{\pi}\biggr)\,+\,\frac{149.9}{\sqrt{\ln{(1/\delta)}}}.
\]
Further, we represent $C_{2}$ as
\[
\frac{1}{\pi-\delta}\sum\limits_{N_{0}<n\le
N}\frac{\Lambda_{1}(n)}{\sqrt{n}}\,n^{-it}\widehat{K}_{a}\biggl(\frac{\ln{n}}{\pi}\biggr)\,-\,\frac{1}{\pi}\sum\limits_{N_{0}<n\le
N}\frac{\Lambda_{1}(n)}{\sqrt{n}}\,n^{-it}\,d_{n},
\]
where
\begin{multline*}
d_{n}\,=\,\widehat{K}_{a}\biggl(\frac{\ln{n}}{\pi}\biggr)\,-\,\widehat{K}_{a}\biggl(\frac{\ln{n}}{\pi-\delta}\biggr)\,=\,
\int_{-\infty}^{+\infty}K_{a}(u)\bigl(e^{-i\varphi_{1}}\,-\,e^{-i\varphi_{2}}\bigr)\,du,\\
\varphi_{1}\,=\,\frac{u\ln{n}}{\pi},\quad
\varphi_{1}\,=\,\frac{u\ln{n}}{\pi-\delta}.
\end{multline*}
Since
\[
\bigl|e^{-i\varphi_{1}}\,-\,e^{-i\varphi_{2}}\bigr|\,=\,2\biggl|\sin{\frac{\varphi_{1}-\varphi_{2}}{2}}\biggr|\,\leqslant\,|\varphi_{1}-\varphi_{2}|\,=\,
\frac{\delta|u|\ln{n}}{\pi(\pi-\delta)},
\]
we obtain:
\[
|d_{n}|\,\leqslant\,\frac{\delta|u|\ln{n}}{\pi(\pi-\delta)}\int_{-\infty}^{+\infty}|u|e^{-\cosh{(\pi
u)}}\,du\,<\,0.01\delta\ln{n}.
\]
Using the bound (\ref{lab_06}) again, we get:
\begin{multline*}
\biggl|\sum\limits_{N_{0}<n\leqslant
N}\frac{\Lambda_{1}(n)}{\sqrt{n}}\,n^{it}\,d_{n}\biggr|\,\leqslant\,0.01\delta\sum\limits_{N_{0}<n\leqslant
N}\frac{\Lambda(n)}{\sqrt{n}}\,\leqslant\\
\leqslant\,0.01\delta\biggl(\frac{\psi(N)}{\sqrt{N}}\,+\,\frac{1}{2}\int_{1}^{N}\frac{\psi(u)}{u^{3/2}}\,du\biggr)\,\leqslant\,0.02c_{1}\delta\sqrt{N}\,<\,
\frac{0.1}{\sqrt{\ln{(1/\delta)}}},
\end{multline*}
and hence
\[
C_{2}\,=\,\frac{1}{\pi-\delta}\sum\limits_{N_{0}<n\le
N}\frac{\Lambda_{1}(n)}{\sqrt{n}}\,n^{-it}\widehat{K}_{a}\biggl(\frac{\ln{n}}{\pi}\biggr)\,+\,
\frac{0.1\theta}{\sqrt{\ln{(1/\delta)}}}.
\]
Finally, the error arising from the replacement of $\pi-\delta$ by
$\pi$ in the last expression does not exceed
\[
\frac{\delta}{\pi(\pi-\delta)}\sum\limits_{N_{0}<n\leqslant
N}\frac{\Lambda_{1}(n)}{\sqrt{n}}\biggl|\widehat{K}_{a}\biggl(\frac{\ln{n}}{\pi}\biggr)\biggr|\,\leqslant\,
\frac{61.5\delta\sqrt{\pi}}{\pi(\pi-\delta)}\sum\limits_{n\leqslant
N}\frac{\Lambda(n)}{n(\ln{n})^{3/2}}\,<\,25\delta
\]
in modulus. Therefore,
\[
C_{2}\,=\,\frac{1}{\pi}\sum\limits_{N_{0}<n\le
N}\frac{\Lambda_{1}(n)}{\sqrt{n}}\,n^{-it}\widehat{K}_{a}\biggl(\frac{\ln{n}}{\pi}\biggr)\,+\,\theta\biggl(25\delta\,+\,
\frac{0.1}{\sqrt{\ln{(1/\delta)}}}\biggr).
\]
Thus, the relation (\ref{lab_05}) takes the form
\begin{multline}\label{lab_07}
\int_{-\infty}^{+\infty}K_{a}\bigl((\pi-\delta)u\bigr)\ln{\zeta\bigl(0.5+i(t+u)\bigr)}\,du\,=\\
=\,\frac{1}{\pi-\delta}\sum\limits_{n \leqslant
N_{0}}\frac{\Lambda_{1}(n)}{\sqrt{n}}\,n^{-it}\widehat{K}_{a}\biggl(\frac{\ln{n}}{\pi-\delta}\biggr)\,+\,
\frac{1}{\pi}\sum\limits_{n>N_{0}}\frac{\Lambda_{1}(n)}{\sqrt{n}}\,n^{-it}\widehat{K}_{a}\biggl(\frac{\ln{n}}{\pi}\biggr)\,-\\
-\,2\pi\int_{0}^{0.5}K_{a}\bigl((\pi-\delta)(t+iv)\bigr)\,dv\,+\,\theta\biggl(25\delta\,+\,
\frac{150}{\sqrt{\ln{(1/\delta)}}}\biggr).
\end{multline}
The integrals in both sides of (\ref{lab_07}) and the sum $C_{3}$
over $n\leqslant N_{0}$ are continuous functions of $\delta$,
$0\leqslant \delta\leqslant 10^{-6}$. Tending $\delta$ to zero, we
lead to the desired statement. The Lemma is proved.

\vspace{0.5cm}

\textbf{\S 2. Basic lemma}

\vspace{0.5cm}

The classical `Dirichlet's approximation theorem' asserts that for
any fixed vector $(\alpha_{1},\ldots,\alpha_{m})$ with real
components and for any arbitrary small $\varepsilon$,
$0<\varepsilon<0.5$, the interval $(1,c)$, $c = \varepsilon^{-m}$,
contains a number $t$ such that the following inequalities hold:
$\|t\alpha_{j}\|<\varepsilon$, $j = 1,\ldots,m$.

Its standard proof (see, for example, \cite[Appendix, \S 9, Theorem
4]{Voronin_Karatsuba_1994}) does not allow one to state the
existence of a number $t$ with the above property on every interval
of the type $(T, T+c_{1})$, where $c_{1}>0$ is a constant depending
only on the tuple $(\alpha_{1},\ldots,\alpha_{m})$ and
$\varepsilon$.

In this section, we prove the analogue of Dirichlet's theorem which
is free of the above disadvantage\footnote{The author sincerely
appreciates O.N.~German and N.G.~Moshchevitin who kindly
communicated him the idea of the proof of Lemma 6.}. However, we
note that the replacement of the interval $(1,c)$ by an arbitrary
interval $(T, T+c_{1})$ leads to the loss of generality (the
condition of linear independence of numbers $1,\alpha_{1},\ldots,
\alpha_{m}$ over the field $\mathbb{Q}$ of the rationals appears)
and to inefficiency of the constant $c_{1} =
c_{1}(\alpha_{1},\ldots, \alpha_{m};\varepsilon)$. The last fact is
a reason of the inefficiency of the constants $c_{0}$ and $T_{0} $
in Theorems 1-6 ($c_{0}$ and $N_{0}$ in Theorem 7,
res\-pec\-ti\-ve\-ly) and of the impossibility of replacement the
value $A$ in Theorem 1 by some increasing function of the parameter
$T$.

\vspace{0.2cm}

\textsc{Lemma 6.} \emph{For any vector $\overline{\alpha} =
(1,\alpha_{1},\ldots,\alpha_{n})$ whose components are linearly
in\-de\-pen\-dent over the rationals and for any $\varepsilon$,
$0<\varepsilon<0.5$, there exists a constant $c =
c(\overline{\alpha},\varepsilon)$ such that any interval of length
$c$ contains at least one value $t$ such that the following
inequalities hold: $\|t\alpha_{j}\|<\varepsilon$, $j = 1,\ldots,
n$.}

\vspace{0.2cm}

\textsc{Proof.} We precede the proof by some remarks.

\textsc{Remark 1.} Let $l$ be the line in $\mathbb{R}^{n+1}$ which
is parallel to vector $\overline{\alpha}$ and passing through the
origin, and let $X = (x_{0}, x_{1},\ldots, x_{n})$ be a point. Then
the distance $d = d(X)$ between $X$ and $l$ is given by a formula
\begin{equation}\label{lab_08}
d\,=\,\frac{1}{|\overline{\alpha}|}\sqrt{\sum\limits_{0\leqslant
i<j\leqslant n}\Delta_{ij}^{2}},\quad\text{where}\quad
|\overline{\alpha}|\,=\,\sqrt{1+\sum\limits_{1\leqslant j\leqslant
n}\alpha_{j}^{2}},
\end{equation}
and $\Delta_{ij}$ is a minor of matrix
\[
\begin{pmatrix}
1 & \alpha_{1} & \dots & \alpha_{n} \\
x_{0} & x_{1} & \dots & x_{n}
\end{pmatrix},
\]
generated by columns $i$ and $j$. Suppose that the lattice point $M
= (m_{0}, m_{1},\ldots, m_{n})$ satisfies the inequality
$d(M)<\varepsilon_{1} = \varepsilon|\alpha|^{-1}$. Then
\[
\sum\limits_{0\leqslant i<j\leqslant
n}\Delta_{ij}^{2}\,<\,\varepsilon^{2}
\]
and therefore
\begin{equation}\label{lab_09}
|\Delta_{01}|\,=\,|\alpha_{1}m_{0}-m_{1}|\,<\,\varepsilon,\quad\ldots,\quad
|\Delta_{0n}|\,=\,|\alpha_{1}m_{0}-m_{n}|\,<\,\varepsilon.
\end{equation}
In view of condition $0<\varepsilon<0.5$, the inequalities
(\ref{lab_09}) imply that $\|\alpha_{j}t\|<\varepsilon$ for any $j$,
$1\leqslant j\leqslant n$, and for $t = m_{0}$.

Thus, it suffices to prove the existence of the infinite sequence of
points $M_{j}$ of the lattice $\mathbb{Z}^{n+1}$ such that the
distance between any neighbouring points $M_{j}$ and $M_{j+1}$ is
bounded from above by some constant depending only on
$\overline{\alpha}$ and $\varepsilon$.

\textsc{Remark 2.} Let us put
\[
\delta\,=\,\frac{\varepsilon_{1}}{n+1}\,=\,\frac{\varepsilon|\overline{\alpha}|^{-1}}{n+1}
\]
and denote by $C_{\delta}$ the infinite cylinder of radius $\delta$
with axis $l$ in $\mathbb{R}^{n+1}$. Suppose that there exist the
points $K_{1},\ldots, K_{n+1} \in \mathbb{Z}^{n+1}$ inside
$C_{\delta}$ such that the vectors $\overline{v}_{j} =
\overrightarrow{OK_{j}\,}$, $j = 1,\ldots, n+1$ are linearly
independent. Then $\overline{v}_{1},\ldots, \overline{v}_{n+1}$
generate an integer lattice $\mathcal{L}$ in $\mathbb{R}^{n+1}$ with
fundamental domain $\Pi$, where $\Pi$ is a parallelepiped spanned on
$\overline{v}_{1},\ldots, \overline{v}_{n+1}$.

It is known that any shift $\Pi + \overline{\xi}$ of the
parallelepiped $\Pi$ to vector $\overline{\xi}\in \mathbb{R}^{n+1}$
contains a point of lattice $\mathcal{L}$ which is also a point of
lattice $\mathbb{Z}^{n+1}$. Further, $\Pi$ is obviously contained in
a cylinder $C_{\varepsilon_{1}} = (n+1)C_{\delta}$ of radius
$(n+1)\delta = \varepsilon_{1}$ which is coaxial to $C_{\delta}$.

Hence, any shift $\Pi + \overline{\xi}$ to vector $\overline{\xi}$
parallel to $\overline{\alpha}$ is fully contained inside
$C_{\varepsilon_{1}}$. At the same time, this shift contains some
lattice point $M(\overline{\xi})$.

Choosing the vectors $\overline{\xi}_{j}$ in such way that the
shifts $\Pi + \overline{\xi}_{j}$ have no pairwise intersections, we
find the desired infinite sequence of lattice points $M_{j} =
M(\overline{\xi}_{j})$ (see Fig. 1).

\begin{center}
\includegraphics{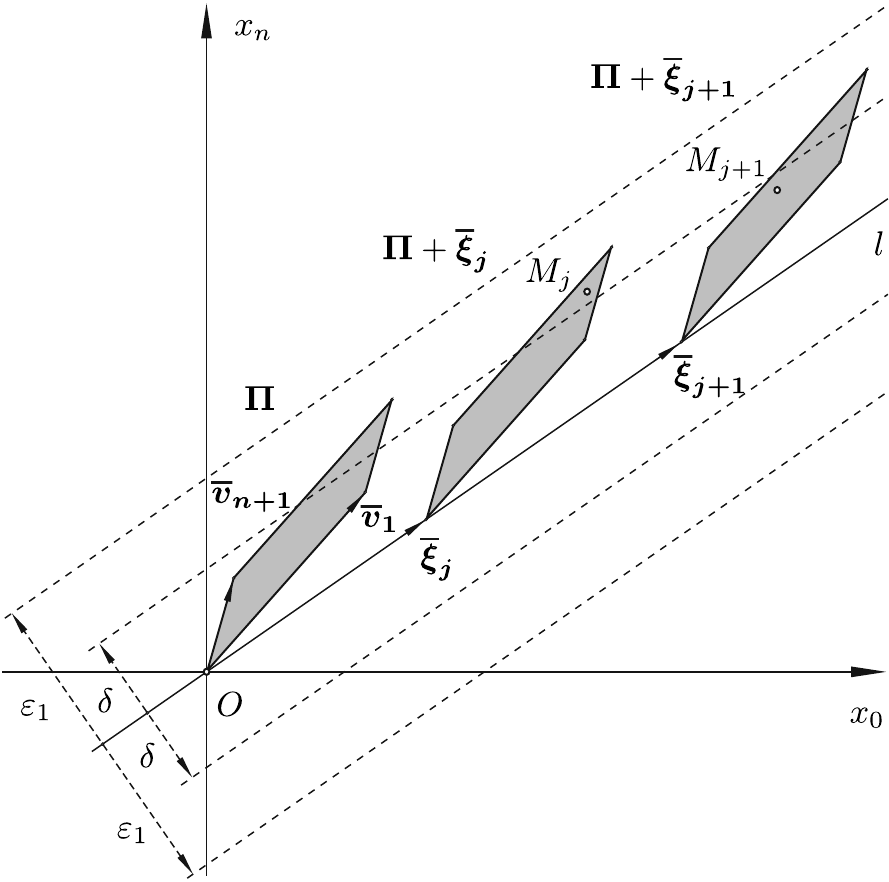}

\fontsize{10}{10pt}\selectfont Fig. 1. Any shift
$\Pi+\overline{\xi}_{j}$ of the parallelepiped $\Pi$ contains a
point $M_{j}$ of the lattice $\mathbb{Z}^{n+1}$.
\fontsize{12}{15pt}\selectfont
\end{center}

Thus, taking $\overline{\xi}_{j} = jc_{0}\,\overline{\alpha}$, $j =
0,\pm 1,\pm 2,\ldots$, where $c_{0} =
2(|\overline{v}_{1}|+\ldots+|\overline{v}_{n+1}|)$ is duplicated sum
of lengths of edges of the parallelepiped $\Pi$ originating from the
same vertex, one can check that the first coordinate of vertex
$\overline{\xi}_{j}$ of $\Pi + \overline{\xi}_{j}$, which is equal
to $jc_{0}$, differs from the first coordinate $m_{0}^{(j)}$ of
lattice point $M_{j}$ for at most
$|\overline{v}_{1}|+\ldots+|\overline{v}_{n+1}| = 0.5c_{0}$. In view
of Remark 1, each of these first coordinates satisfies the series of
inequalities $\|\alpha_{i}m_{0}^{(j)}\|<\varepsilon$, $i = 1,\ldots,
n+1$. Since
\[
\bigl|m_{0}^{(j)}\,-\,m_{0}^{(j+1)}\bigr|\,\leqslant\,(j+1)c_{0}+0.5c_{0}-(jc_{0}-0.5c_{0})\,=\,2c_{0},
\]
it appears that any interval of the type $(\tau, \tau+3c_{0})$
contains a point of sequence $m_{0}^{(j)}$, $j = 0,\pm 1,\pm
2,\ldots$.

Thus, it suffices to prove that any cylinder $C_{\delta}$ with axis
$l$ contains $n+1$ linearly independent vectors of the lattice
$\mathbb{Z}^{n+1}$.

Now let us prove the main assertion. First we show that $C_{\delta}$
contains an infinite set of lattice points.

The line $l$ does not contain lattice points different from the
origin $O$. In the opposite case, we have $d(K) = 0$, $k_{0}\ne 0$
for such point $K = (k_{0}, k_{1},\ldots, k_{n})\in
\mathbb{Z}^{n+1}$. Hence $\Delta_{0j}=\alpha_{j}k_{0} - k_{j} = 0$
for any $j = 1,\ldots, n$ and therefore, $\alpha_{j} =
k_{j}/k_{0}\in \mathbb{Q}$. But this contradicts the linear
independence of $1,\alpha_{1}, \ldots, \alpha_{n}$ over the
rationals.

Let $\Omega_{n}$ be an $n$\,-dimensional hyperplane passing through
the origin $O$ perpendicularly to the axis $l$. Then an
$n$\,-dimensional volume $V_{1}$ of a sphere arising in the intersection of
the cylinder $C_{\delta}$ with the hyperplane $\Omega_{n}$ is equal to
$V_{1} = c(n)\delta^{\,n}$, where $c(n) =
\pi^{n/2}\Gamma^{\,-1}\bigl(n/2+1\bigr)$. Now let us define $H_{1}$
by the relation $H_{1}V_{1} = 2^{n-1}$ and consider an
$(n+1)$\,-dimensional cylinder $T_{1}$ of height $2H_{1}$ which
arises from $C_{\delta}$ after the section by two hyperplanes parallel to
$\Omega_{n}$ which are distant to $H_{1}$ from the origin.

Since the volume of such cylinder is equal to $2H_{1}V_{1} = 2^{n}$,
Minkowski's convex body theorem (see for example \cite[\S
5]{Gruber_2008}) implies that this cylinder contains a lattice point
$N_{1}$ different from the origin $O$.

\begin{center}
\includegraphics{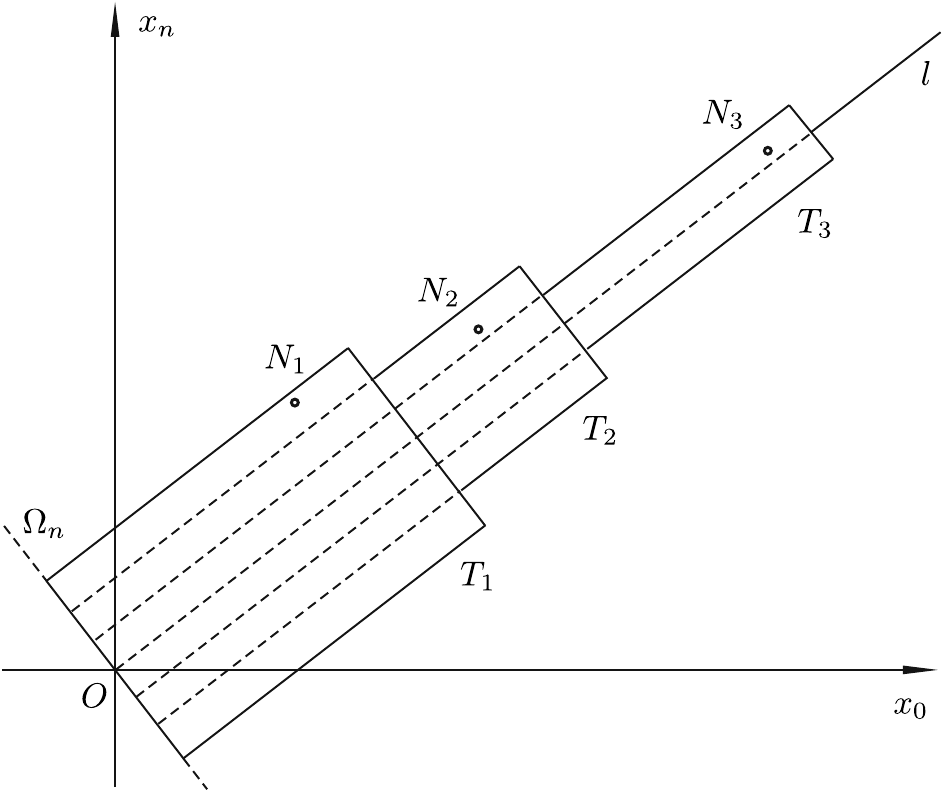}

\fontsize{10}{10pt}\selectfont Fig. 2. An infinite sequence of
lattice points $N_{j}$. \fontsize{12}{15pt}\selectfont
\end{center}

Without loss of generality, we assume that $N_{1}$ is the closest
point to $l$ among the lattice points of the cylinder $T_{1}$ which
differs from the origin $O$. In view of the above remark, $N_{1}$
does not lie on $l$, so we have $d(N_{1})>0$.

Further, let us take $\delta_{2} = 0.5d(N_{1})$ and define $H_{2}$
by the relations $H_{2}V_{2} = 2^{n-1}$, $V_{2} =
c(n)\delta_{2}^{n}$. Applying the same arguments to the cylinder
$T_{2}$ of radius $\delta_{2}$ and height $2H_{2}$, which is
symmetrical with respect to the origin and coaxial to $T_{1}$, we
find a lattice point $N_{2}$ inside it, which is different from the
origin $O$ and closest to $l$ among the lattice points of $T_{2}$.
Since $d(N_{2})\leqslant\delta_{2}<d(N_{1})$, the point $N_{2}$
differs from $N_{1}$. In view of symmetry both of $T_{1}$ and
$T_{2}$ with respect to $O$, we assume that $N_{1}$ and $N_{2}$ lie
in the same half\,-space with respect to the hyperplane
$\Omega_{n}$.

Taking $\delta_{3} = 0.5d(N_{2})$, $H_{3}V_{3} = 2^{n-1}$, $V_{3} =
c(n)\delta_{3}^{\,n}$, we construct in the same way the cylinder
$T_{3}$ of radius $\delta_{3}$ and height $2H_{3}$ and find a
lattice point $N_{3}$ inside it, which differs from $O$, $N_{1}$,
$N_{2}$ and lying in the same half\,-space with respect to
$\Omega_{n}$.

Proceeding this process further, we finally get an infinite sequence
of different points $N_{j}$ of the lattice $\mathbb{Z}^{n+1}$
containing in the same half of the cylinder $C_{\delta}$ with
respect to secant hyperplane $\Omega_{n}$ and satisfying the
condition $0<d(N_{j+1})\leqslant 0.5d(N_{j})$, $j = 1,2,3,\ldots$.

Now we prove the existence of $n+1$ linearly independent vectors among the infinite set $\overrightarrow{ON_{j}\,}$, $j = 1,2,3,\ldots$.

Let's assume the contrary. Suppose that the maximal number $s$ of linearly independent
vectors from this set does not exceed $n$. Let $\overline{u}_{1}, \ldots , \overline{u}_{s}\in \mathbb{Z}^{n+1}$ be such vectors
and let $\omega_{s}$ be the $s$\,-dimensional hyperplane spanned on it.

Then the intersection of $\omega_{s}$ and $C_{\delta}$ contains an
infinite sequence of points $N_{j}$ of lattice $\mathbb{Z}^{n+1}$.
Hence, this intersection is unbounded. But the intersection of
$\omega_{s}$ and $C_{\delta}$ is unbounded if and only if the
hyperplane $\omega_{s}$ is parallel to the line $l$ or contains it
(see Fig. 3).

\begin{center}
\includegraphics{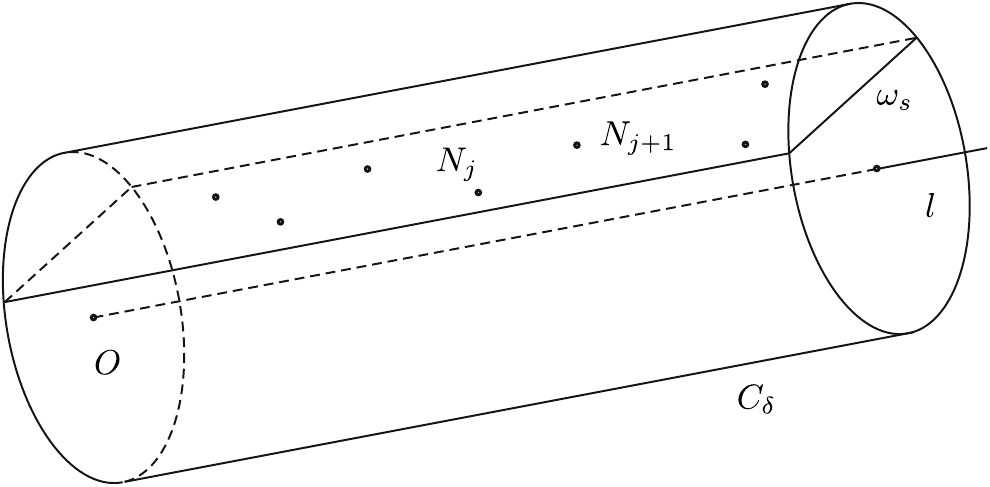}

\fontsize{10}{10pt}\selectfont Fig. 3. The intersection of
$C_{\delta}$ and $\omega_{s}$ is unbounded.
\fontsize{12}{15pt}\selectfont
\end{center}

If the first case, all the distances between $N_{j}$ and $l$ are
bounded from below by some positive constant (which is equal to the distance between $\omega_{s}$ and $l$).
But this is impossible since $d(N_{j})\to 0$ as $j\to +\infty$.

Further, if the line $l$ lies in the hyperplane $\omega_{s}$ then $\overline{\alpha}$ is the linear combination of the form
$\overline{\alpha}=t_{1}\overline{u}_{1}+\ldots+t_{s}\overline{u}_{s}$. Denoting the components of $\overline{u}_{j}$ by
$u_{0j}, u_{1j},\ldots, u_{nj}$, we get:
\begin{equation}\label{lab_10}
\left\{
\begin{aligned}
t_{1}u_{01}\,+\,\ldots\,+\,t_{s}u_{0s}\,& =\,1,\\
t_{1}u_{11}\,+\,\ldots\,+\,t_{s}u_{1s}\,& =\,\alpha_{1},\\
\cdots \\
t_{1}u_{n1}\,+\,\ldots\,+\,t_{s}u_{ns}\,& =\,\alpha_{n}.
\end{aligned}
\right.
\end{equation}
Since $\overline{u}_{1},\ldots, \overline{u}_{s}$ are linearly
independent then $(n+1)\times s$\,-matrix of its components has the
maximal rank $s$. Hence, it contains $s$ linearly independent rows,
and let $0\leqslant i_{1}<i_{2}<\ldots < i_{s}\leqslant n$ be their
indices. If it is necessary, we put $\alpha_{0} = 1$ and consider
the corresponding system of equations extracting from
(\ref{lab_10}), that is
\[
\left\{
\begin{aligned}
t_{1}u_{i_{1}1}\,+\,\ldots\,+\,t_{s}u_{i_{1}s}\,& =\,\alpha_{i_{1}},\\
\cdots \\
t_{1}u_{i_{s}1}\,+\,\ldots\,+\,t_{s}u_{i_{s}s}\,& =\,\alpha_{i_{s}}.
\end{aligned}
\right.
\]
Its determinant is nonzero integer. Cramer's formulas implies that the unique solution of this system has the form
\[
\left\{
\begin{aligned}
t_{1}\, &=\,r_{11}\alpha_{i_{1}}\,+\,\ldots\,+\,r_{1s}\alpha_{i_{s}},\\
\cdots \\
t_{s}\,
&=\,r_{s1}\alpha_{i_{1}}\,+\,\ldots\,+\,r_{ss}\alpha_{i_{s}},
\end{aligned}
\right.
\]
where $r_{ij}$ are some rationals. Since $s\leqslant n$ then there
exist at least one equation in (\ref{lab_10}) whose index $j$
differs from $i_{1},\ldots, i_{s}$. Thus we get:
\begin{multline*}
\alpha_{j}\,=\,t_{1}u_{j1}\,+\,\ldots\,+\,t_{s}u_{js}\,=\\
=\, \bigl(\,r_{11}\alpha_{i_{1}}\,+\,\ldots\,+\,r_{1s}\alpha_{i_{s}}\bigr)u_{j1}\,+\,\ldots\,+\,
\bigl(r_{s1}\alpha_{i_{1}}\,+\,\ldots\,+\,r_{ss}\alpha_{i_{s}}\bigr)u_{js}\,=\\
=\,q_{1}\alpha_{i_{1}}\,+\,\ldots\,+\,q_{s}\alpha_{s},
\end{multline*}
where $q_{1},\ldots, q_{s}\in \mathbb{Q}$. The last relation contradicts to the linear independence
of $1,\alpha_{1},\ldots,$ $\alpha_{n}$ over the rationals.

This contradiction implies that the hyperplane $\omega_{s}$ does not contain the line $l$. This proves the lemma.

\vspace{0.2cm}

\textsc{Corollary.} \emph{For any vector $\overline{\alpha} =
(1,\alpha_{1},\ldots,\alpha_{n})$ whose components are linearly
independent over the rationals, for any tuple of real numbers
$\beta_{1},\ldots,$ $\beta_{n}$ and for any $\varepsilon$,
$0<\varepsilon<0.5$, there exists a constant $c =
c(\overline{\alpha},\varepsilon)$ such that any interval of length
$c$ contains at least one value $t$ such that the following
inequalities hold: $\|t\alpha_{j}+\beta_{j}\|<\varepsilon$, $j =
1,\ldots, n$.}

\vspace{0.2cm}

\textsc{Proof.} We use the notations of Lemma 6. The above arguments
imply that the cylinder $C$ with radius $\varepsilon_{1} =
\varepsilon|\overline{\alpha}|^{-1}$ and axis $l$ passing through
the origin in parallel to $\overline{\alpha}$ contains an
$(n+1)$\,-dimensional parallelepiped $\Pi$ whose vertices belong to
the lattice $\mathbb{Z}^{n+1}$.

Then the cylinder $C_{0} = C+\overline{\beta}$, which is the shift
of $C$ to vector $\overline{\beta} =
(1,\beta_{1},\ldots,\beta_{n})$, contains a parallelepiped $\Pi_{0}
= \Pi + \overline{\beta}$. Any shift of $\Pi$ contains a lattice
point. Hence, both $\Pi_{0}$ and any parallelepiped $\Pi_{j}$ which
is the shift of $\Pi_{0}$ to vector $\overline{\xi}_{j} =
c_{0}j\,\overline{\alpha}$, $j = \pm 1, \pm 2, \ldots$, parallel to
the axis of the cylinder $C_{0}$, contain the points of the lattice
$\mathbb{Z}^{n+1}$. It is easy to note that the parallelepipeds
$\Pi_{j}$ have no common points.

Finally, let $M_{j} = (m_{0},\ldots, m_{n})$  be a lattice point
containing in $\Pi_{j}$. The distance between this point and the
axis of $C_{0}$ does not exceed $\varepsilon_{1}$. At the same time,
this distance is expressed by (\ref{lab_08}), where $\Delta_{ij}$ is
a minor of matrix
\[
\begin{pmatrix}
1 & \alpha_{1} & \dots & \alpha_{n} \\
m_{0} & m_{1}-\beta_{1} & \dots & m_{n}-\beta_{n}
\end{pmatrix}.
\]
formed by its columns $i, j$. Hence, we have
\[
|\Delta_{ij}|\,=\,\bigl|m_{0}\alpha_{j}\,-\,(m_{j}-\beta_{j})\bigr|\,=\,\bigl|m_{0}\alpha_{j}+\beta_{j}-m_{j}\bigr|\,<\,\varepsilon
\]
for any $j$, $1\leqslant j\leqslant n$. By the inequality
$\varepsilon<0.5$, we obtain that
$\|m_{0}\alpha_{j}+\beta_{j}\|<\varepsilon$. To end the proof, we
note that the first coordinates $m_{0}$ of the points $M_{j}$ form
an increasing sequence, whose neighbouring elements differ for at
most to $3c_{0}$.

\vspace{0.5cm}

\textbf{\S 3. Large values of the Riemann zeta function on the
critical line}

\vspace{0.5cm}

In this section, we give a conditional solution of Karatsuba's problem based on the Riemann hypothesis.
We also prove a series of statements concerning the existence of large values of the function $S(t)$ on the short
segments of the real axis.

\vspace{0.2cm}

\textsc{Theorem 1.} \emph{Suppose that the Riemann hypothesis is
true, and let $A$ be an arbitrary large fixed constant. Then there
exist the constants $c_{0} = c_{0}(A)>0$ and $T_{0} = T_{0}(A)$ such
that any interval of the form $(T-H,T+H)$, $H =
(1/\pi)\ln\ln\ln{T}+c$, $T>T_{0}$, contains at least one point $t$
such that $\bigl|\zeta(0.5+it)\bigr|>A$.}

\vspace{0.2cm}

\textsc{Proof.} Let's fix any positive number
$a>1$ satisfying the condition
\begin{equation}\label{lab_11}
e^{a}\sqrt{\frac{\pi}{2a}}\,\geqslant\,\ln{A}.
\end{equation}
Extracting real parts in (\ref{lab_04}), we obtain:
\begin{multline}\label{lab_12}
\int_{-\infty}^{+\infty}K_{a}(\pi
u)\ln{\bigl|\zeta\bigl(0.5+i(t+u)\bigr)\bigr|}\,du\,=\,\frac{1}{\pi}\sum\limits_{n
=
2}^{+\infty}\frac{\Lambda_{1}(n)}{\sqrt{n}}\,\widehat{K}_{a}\biggl(\frac{\ln{n}}{\pi}\biggr)\cos{(t\ln{n})}\,-\\
-\,2\pi\int_{0}^{0.5}\Re{K_{a}(\pi t + \pi i v)}dv.
\end{multline}
Taking $t = 0$ in (\ref{lab_12}) and noting that $K_{a}(\pi
iv)=e^{-a\cos{\pi v}}$, we have:
\begin{equation}\label{lab_13}
\int_{-\infty}^{+\infty}K_{a}(\pi
u)\ln{\bigl|\zeta\bigl(0.5+iu\bigr)\bigr|}\,du\,=\,\frac{1}{\pi}\sum\limits_{n
=
2}^{+\infty}\frac{\Lambda_{1}(n)}{\sqrt{n}}\,\widehat{K}_{a}\biggl(\frac{\ln{n}}{\pi}\biggr)\,-\,2\pi\int_{0}^{0.5}e^{-a\cos{(\pi
v)}}\,dv.
\end{equation}
Further, the relation $\bigl|K_{a}(\pi t+\pi i
v)\bigr|=e^{-a\cosh{(\pi t)\cos{(\pi v)}}}$ implies that the last
integral in (\ref{lab_12}) does not exceed
\begin{equation}\label{lab_14}
2\pi\int_{0}^{0.5}e^{-a\cosh{(\pi t)}\cos{(\pi
v)}}\,dv\,=\,2\pi\int_{0}^{0.5}e^{-a\cosh{(\pi t)}\sin{(\pi
v)}}\,dv\,<\,\frac{\pi}{a\cosh{(\pi t)}}
\end{equation}
in modulus. Subtracting (\ref{lab_13}) from (\ref{lab_12}) and using the estimate
(\ref{lab_14}), we find
\begin{multline}\label{lab_15}
\int_{-\infty}^{+\infty}K_{a}(\pi
u)\ln{\bigl|\zeta\bigl(0.5+i(t+u)\bigr)\bigr|}\,du\,-\,\int_{-\infty}^{+\infty}K_{a}(\pi
u)\ln{\bigl|\zeta\bigl(0.5+iu\bigr)\bigr|}\,du\,=\\
=\,2\pi\int_{0}^{0.5}e^{-a\cos{(\pi
v)}}\,dv\,-\,\frac{2}{\pi}\sum\limits_{n =
2}^{+\infty}\frac{\Lambda_{1}(n)}{\sqrt{n}}\,\widehat{K}_{a}\biggl(\frac{\ln{n}}{\pi}\biggr)\sin^{2}\biggl(\frac{t}{2}\ln{n}\biggr)\,+\,
\frac{\pi\theta_{1}}{\cosh{(\pi t)}}.
\end{multline}

Let $\varepsilon, N$ be the numbers satisfying the conditions
$0<\varepsilon<0.5$, $N>N_{0}=e^{\pi a\sqrt{2}}$ and depending only
on $a$, whose precise values will be chosen below. Applying Lemmas 1
and 6, we find the constant  $c_{0} = c_{0}(a)$ such that any
interval of real axis with length $c_{0}$ contains at least one
point $\tau$ such that the inequalities
$\bigl\|(\tau/(2\pi))\ln{p}\bigr\|<\varepsilon$ hold true for all
primes $p\leqslant N$. Let us take $t$ to be equal to such value
$\tau$ from the interval $(T,T+c_{0})$ in (\ref{lab_15}).

Given prime $p\leqslant N$, we define an integer $n_{p}$ and real
$\varepsilon_{p}$ satisfying the condition
$|\varepsilon_{p}|<\varepsilon$ such that $(t/(2\pi))\ln{p} =
n_{p}+\varepsilon_{p}$. Then we have
\[
\sin^{2}\biggl(\frac{t}{2}\ln{n}\biggr)\,=\,\sin^{2}\bigl(\pi
kn_{p}+\pi k\varepsilon_{p}\bigr)\,=\,\sin^{2}(\pi
k\varepsilon_{p})\,<\,(\pi k\varepsilon)^{2}
\]
for any $k\geqslant 1$ and $n = p^{k}$.

Let $C$ be the sum in the right\,-hand side of (\ref{lab_15}).
Denote by $C_{1}$ and $C_{2}$ the contributions to $C$ from the
terms corresponding to $n = p^{k}$, $k\geqslant 1$, $p\leqslant N$
and from all other terms, respectively. Then we have:
\[
|C_{1}|\,\leqslant\,\frac{2}{\pi}\,(\pi\varepsilon)^{2}\!\!\sum\limits_{\substack{n\,=\,
p^{k} \\ k\geqslant 1, \,p\leqslant
N}}\frac{k}{\sqrt{n}}\biggl|\widehat{K}_{a}\biggl(\frac{\ln{n}}{\pi}\biggr)\biggr|.
\]
We split the domain of $n$ to the intervals $n\leqslant N_{0}$,
$N_{0}<n\leqslant N$ and $n>N$ and then denote the corresponding
parts of sum by $C_{3}, C_{4}, C_{5}$. The estimate
$\bigl|\widehat{K}_{a}\bigl((1/\pi)\ln{n}\bigr)\bigr|\,\leqslant\,\widehat{K}_{a}(0)$
implies
\begin{multline*}
|C_{3}|\,\leqslant\,2\pi\varepsilon^{2}\widehat{K}_{a}(0)
\sum\limits_{p\leqslant N_{0}}\sum\limits_{k = 1}^{+\infty}kp^{-\,k/2}\,=\\
=\,2\pi\varepsilon^{2}\widehat{K}_{a}(0)\sum\limits_{p\leqslant
N_{0}}\frac{1}{\sqrt{p}}\biggl(1\,-\,\frac{1}{\sqrt{p}}\biggr)^{-2}\,\leqslant\,2\pi\varepsilon^{2}
\biggl(1-\frac{1}{\sqrt{2}}\biggr)^{-2}\widehat{K}_{a}(0)\sum\limits_{p\leqslant
N_{0}}\frac{1}{\sqrt{p}}.
\end{multline*}
Let us use the inequality
\[
\sum\limits_{p\,\leqslant\,x}\frac{1}{\sqrt{p}}\,\leqslant\,\frac{2.784\sqrt{x}}{\ln{x}},
\]
which is verified for $2\leqslant x\leqslant 1.5\cdot 10^{6}$ by
Wolfram Mathematica 7.0 and follows from the inequality (3.6) from
\cite[Th. 2, corollary 1]{Rosser_Schoenfeld_1962} by Abel's
summation formula for $x>1.5\cdot 10^{6}$. Thus we get
\[
|C_{3}|\,<\,45.9\varepsilon^{2}\widehat{K}_{a}(0)\,\frac{e^{\pi
a/\sqrt{2}}}{a}\,\leqslant\,(7\varepsilon)^{2}e^{\pi
a/\sqrt{2}}\widehat{K}_{a}(0).
\]
Further, the Corollary of Lemma 3 implies
\begin{multline*}
|C_{4}|\,\leqslant\,2\pi\varepsilon^{2}\sum\limits_{\substack{N_{0}<n\leqslant N \\
n\,=\,p^{k}}}
kp^{-k/2}\cdot\frac{61.5\sqrt{\pi}}{\sqrt{k\ln{p}\mathstrut}}\,
\exp{\biggl(-\frac{\pi}{2}\,\frac{1}{\pi}\ln{p^{k}}\biggr)}\,=\\
=\,123\pi\sqrt{\pi}\varepsilon^{2}\sum\limits_{\substack{N_{0}<n\leqslant N \\
n\,=\,p^{k}}}
\frac{\sqrt{k}}{p^{k\mathstrut}\sqrt{\ln{p}\mathstrut}}\,\leqslant\,
123\pi\sqrt{\pi}\varepsilon^{2}\sum\limits_{p\leqslant N}\frac{1}{\sqrt{\ln{p}\mathstrut}}\sum\limits_{k = 1}^{+\infty}kp^{-k}\,<\\
<\,123\pi\sqrt{\pi}\varepsilon^{2}\sum\limits_{p\leqslant
N}\frac{1}{p\sqrt{\ln{p}\mathstrut}}
\biggl(1\,-\,\frac{1}{p}\biggr)^{\!-2}\,<\,123\pi\sqrt{\pi}\varepsilon^{2}\sum\limits_{p}\frac{p}{(p-1)^{2}\sqrt{\ln{p}\mathstrut}}\,<\,3000\varepsilon^{2}.
\end{multline*}
Applying the Corollary of lemma 3 together with the estimate
(\ref{lab_06}) again and noting that $\ln{N}\geqslant$ $\pi
a\sqrt{2}\geqslant\pi\sqrt{2}$, we find
\[
|C_{5}|\,\leqslant\,\frac{2}{\pi}\sum\limits_{n>N}\frac{\Lambda_{1}(n)}{\sqrt{n}}\,
\frac{61.5\sqrt{\pi}}{\sqrt{n\ln{n}\mathstrut}}\,=\,\frac{123}{\sqrt{\pi}}
\sum\limits_{n>N}\frac{\Lambda(n)}{n(\ln{n})^{3/2}}.
\]
Abel's summation formula together with the bound
\[
\psi(u)\,=\,\sum\limits_{n\leqslant u}\Lambda(n)\,\leqslant\,c_{1}u,
\quad c_{1}\,=\,1.03883
\]
(see \cite[Th. 12]{Rosser_Schoenfeld_1962}), which is valid for any
$u>0$, imply
\begin{multline*}
\sum\limits_{n>N}\frac{\Lambda(n)}{n(\ln{n})^{3/2}}\,=\,-\int_{N}^{+\infty}\bigl(\psi(u)\,-\,\psi(N)\bigr)\,
d\,\frac{1}{(\ln{u})^{3/2}}\,\leqslant\,-c_{1}\int_{N}^{+\infty}u\,d\,\frac{1}{(\ln{u})^{3/2}}\,=\\
=\,c_{1}\biggl(\frac{2}{\sqrt{\ln{N}\mathstrut}}\,+\,\frac{1}{(\ln{N})^{3/2}}\biggr).
\end{multline*}
Since $\ln{N}\geqslant\pi a\sqrt{2}\geqslant\pi\sqrt{2}$, we finally
get:
\begin{multline*}
|C_{5}|\,\leqslant\,\frac{123}{\sqrt{\pi}}\,\frac{2c_{1}}{\sqrt{\ln{N}\mathstrut}}\biggl(1\,+\,\frac{1}{2\pi\sqrt{2}}\biggr)\,<\,
\frac{160.5}{\sqrt{\ln{N}\mathstrut}},\\
|C_{1}|\,\leqslant\,|C_{3}|\,+\,|C_{4}|\,+\,|C_{5}|\,<\,(7\varepsilon)^{2}\widehat{K}_{a}(0)e^{\pi
a/\sqrt{2}}\,+\,3000\varepsilon^{2}\,+\,\frac{160.5}{\sqrt{\ln{N}\mathstrut}}.
\end{multline*}
Applying the same arguments to the estimation of the sum $C_{2}$, we obtain
\[
|C_{2}|\,\leqslant\,\frac{2}{\pi}\sum\limits_{n>N}\frac{\Lambda_{1}(n)}{\sqrt{n}}\,\frac{61.5\sqrt{\pi}}{\sqrt{n\ln{n}\mathstrut}}\,<\,
\frac{160.5}{\sqrt{\ln{N}\mathstrut}}.
\]
Thus
\[
|C|\,\leqslant\,|C_{1}|\,+\,|C_{2}|\,<\,(7\varepsilon)^{2}\widehat{K}_{a}(0)e^{\pi
a/\sqrt{2}}\,+\,3000\varepsilon^{2}\,+\,\frac{321}{\sqrt{\ln{N}\mathstrut}},
\]
and therefore
\begin{multline}\label{lab_16}
\int_{-\infty}^{+\infty}K_{a}(\pi u)\ln{\bigl|\zeta\bigl(0.5+i(t+u)\bigr)\bigr|}\,du\,\geqslant\,2\int_{0}^{\pi/2}e^{-\,a\sin{v}}\,dv\,+\\
+\,2\int_{0}^{+\infty}K_{a}(\pi
u)\ln{\bigl|\zeta\bigl(0.5+iu\bigr)\bigr|}\,du\,-\,\biggl((7\varepsilon)^{2}\widehat{K}_{a}(0)e^{\pi
a/\sqrt{2}}\,+\,3000\varepsilon^{2}\,+\,\frac{321}{\sqrt{\ln{N}\mathstrut}}\,+\,\frac{\pi}{\cosh{\pi
t}}\biggr).
\end{multline}

Now we estimate the modulus of the improper integral in the
right\,-hand side of (\ref{lab_16}). We split it to the integrals
$j_{1}$ and $j_{2}$, corresponding  to the intervals $0\leqslant
u\leqslant 10$ and $u>10$, respectively. Since the modulus of
$\ln{\bigl|\zeta(0.5+iu)\bigr|}$  does not exceed $0.641973\ldots <
2/3-1/50$ for $0\leqslant u\leqslant 10$, we find
\[
|j_{1}|\,<\,\biggl(\frac{2}{3}\,-\,\frac{1}{50}\biggr)\int_{0}^{10}K_{a}(\pi
u)\,du\,<\,\frac{1}{\pi}\biggl(\frac{1}{3}\,-\,\frac{1}{100}\biggr)\widehat{K_{a}}(0).
\]
Further, the formula for $\widehat{K}_{a}(0)$ from \cite[Ex. 9.1]{Olver_1986} implies that
\begin{equation}\label{lab_17}
\frac{7}{8}e^{-a}\sqrt{\frac{2\pi}{a}}\,<\,\widehat{K}_{a}(0)\,<\,e^{-a}\sqrt{\frac{2\pi}{a}}
\end{equation}
for $a>1$. Hence,
\begin{multline*}
|j_{2}|\,\leqslant\,\frac{\widehat{K}_{a}(0)}{\widehat{K}_{a}(0)}\int_{10}^{+\infty}e^{-a\cosh{(\pi u)}}\bigl|\ln{\bigl|\zeta(0.5+iu)\bigr|}\bigr|\,du\,\leqslant\\
\leqslant\,\widehat{K}_{a}(0)\,\frac{8}{7}\,e^{a}\sqrt{\frac{a}{2\pi}}\int_{10}^{+\infty}
\exp{\bigl(-\,0.5ae^{\pi u}\bigr)}\bigl|\ln{\bigl|\zeta(0.5+iu)\bigr|}\bigr|\,du\,=\\
=\,\widehat{K}_{a}(0)\,\frac{8}{7}\,e^{-a}\sqrt{\frac{a}{2\pi}}\int_{10}^{+\infty}\exp{\bigl(-\,0.5a\bigl(e^{\pi
u}-4\bigr)\bigr)}\bigl|\ln{\bigl|\zeta(0.5+iu)\bigr|}\bigr|\,du.
\end{multline*}
Since $0.5\bigl(e^{\pi u}-4\bigr)>2u^{2}$ for $u\geqslant 10$, we
find
\begin{multline*}
|j_{2}|\,\leqslant\,\frac{\widehat{K}_{a}(0)}{\widehat{K}_{a}(0)}\int_{10}^{+\infty}e^{-2u^{2}}\bigl|\ln{\bigl|\zeta(0.5+iu)\bigr|}\bigr|\,du\,
\leqslant\,\widehat{K}_{a}(0)\,\frac{8}{7}\,e^{-\,a}\sqrt{\frac{a}{2\pi}}\cdot
1.52\cdot 10^{-89}\,<\\
<\,1.5\cdot 10^{-90}\widehat{K}_{a}(0).
\end{multline*}
Thus we get
\[
|j_{1}|\,+\,|j_{2}|\,<\,\frac{1}{\pi}\biggl(\frac{1}{3}\,-\,\frac{1}{100}\biggr)\widehat{K}_{a}(0)\,+\,1.5\cdot
10^{-90}\,\widehat{K}_{a}(0)\,<\,\frac{\widehat{K}_{a}(0)}{3\pi}.
\]
Obviously we have
\[
\int_{0}^{\pi/2}e^{-a\sin{v}}\,dv\,\geqslant\,\int_{0}^{\pi/2}e^{-av}\,dv\,=\,\frac{1}{a}\bigl(1\,-\,e^{-\pi
a/2}\bigr).
\]
Therefore, the inequality (\ref{lab_16}) implies
\begin{multline}\label{lab_18}
\int_{-\infty}^{+\infty}K_{a}(\pi u)\ln{\bigl|\zeta\bigl(0.5+i(t+u)\bigr)\bigr|}\,du\,\geqslant\,\frac{2}{a}\bigl(1\,-\,e^{-\pi a/2}\bigr)\,-\\
-\,\biggl((7\varepsilon)^{2}\widehat{K}_{a}(0)e^{\pi
a/\sqrt{2}}\,+\,3000\varepsilon^{2}\,+\,\frac{321}{\sqrt{\ln{N}\mathstrut}}\,+\,\frac{\widehat{K}_{a}(0)}{3\pi}\,+\,\frac{\pi}{\cosh{\pi
t}}\biggr).
\end{multline}
Further, we put $h = (1/\pi)(\ln\ln\ln{T}-\ln{(a/2)})$ and split the
integral to the sum
\[
\biggl(\,\int_{-h}^{h}\,+\,\int_{h}^{+\infty}\,+\,\int_{-\infty}^{-h}\,\biggr)K_{a}(\pi
u)
\ln{\bigl|\zeta\bigl(0.5+i(t+u)\bigr)\bigr|}\,du\,=\,j_{3}\,+\,j_{4}\,+\,j_{5}.
\]
The formula
\[
\zeta(s)\,=\,\frac{1}{s-1}\,+\,\frac{1}{2}\,+\,s\int_{1}^{+\infty}\frac{\varrho(u)}{u^{s+1}}\,du,
\]
where $\varrho(u) = 0.5-\{u\}$, $\Re s >0$, $s\ne 1$ (see \cite[Ch.
II, Lemma 2]{Karatsuba_1983}) implies that $0\leqslant
\bigl|\zeta(0.5+iv)\bigr|\leqslant |v|+3$ for any real $v$. Hence,
\[
-\infty\,\leqslant\,\ln{\bigl|\zeta(0.5+iv)\bigr|}\,<\,\ln{\bigl(|v|+3\bigr)}.
\]
Passing to the estimate of $j_{4}$, we get:
\begin{multline*}
-\infty\,\leqslant\,j_{4}\,=\,\int_{h}^{+\infty}K_{a}(\pi
u)\ln{\bigl|\zeta\bigl(0.5+i(t+u)\bigr)\bigr|}\,du\,<\,
\int_{h}^{+\infty}K_{a}(\pi u)\ln{\bigl(|t+u|+3\bigr)}\,du\,=\\
=\,\biggl(\,\int_{h}^{t}\,+\,\int_{t}^{+\infty}\,\biggr)K_{a}(\pi
u)\ln{\bigl(|t+u|+3\bigr)}\,du\,=\,j_{6}\,+\,j_{7}.
\end{multline*}
Estimating the integrals $j_{6}$ и $j_{7}$ separately, we find
\begin{multline*}
j_{6}\,\leqslant\,\ln{(2t+3)}\int_{h}^{+\infty}\exp{\bigl(-0.5ae^{\pi
u}\bigr)}\,du\,=\,\frac{1}{\pi}\ln{(2t+3)}
\int_{0.5ae^{\pi h}}^{+\infty}e^{-w}\frac{dw}{w}\,=\\
=\,\frac{1}{\pi}\ln{(2t+3)}\int_{\ln\ln{T}}^{+\infty}e^{-w}\frac{dw}{w}\,<\,\frac{\ln{(2t+3)}}{\pi\ln{T}}\,\frac{1}{\ln\ln{T}}.
\end{multline*}
Similarly,
\begin{multline*}
j_{7}\,\leqslant\,\int_{t}^{+\infty}\exp{\bigl(-0.5ae^{\pi u}\bigr)}\ln{(2u+3)}\,du\,\leqslant\,2\int_{t}^{+\infty}\exp{\bigl(-0.5ae^{\pi u}\bigr)}(\ln{u})\,du\,<\\
<\,\frac{2}{\pi}\ln{\bigl(\pi t/2\bigr)}e^{-\pi
t/2}\exp{\bigl(-\,e^{\pi t/2}\bigr)}.
\end{multline*}
Therefore,
\[
-\infty\,\leqslant\,j_{4}\,=\,j_{6}\,+\,j_{7}\,<\,
\frac{\ln{(2t+3)}}{\pi\ln{T}}\,\frac{1}{\ln\ln{T}}\,+\,\frac{2}{\pi}\ln{\bigl(\pi
t/2\bigr)}e^{-\pi t/2}\exp{\bigl(-\,e^{\pi
t/2}\bigr)}\,<\,\frac{1}{3\ln\ln{T}}.
\]
The integral $j_{5}$ is estimated in the same way. Thus we have:
\begin{multline*}
j_{5}\,=\,\int_{h}^{+\infty}K_{a}(\pi
u)\ln{\bigl|\zeta\bigl(0.5+i(t-u)\bigr)\bigr|}\,du\,<\,
\int_{h}^{+\infty}K_{a}(\pi u)\ln{\bigl(|t-u|+3\bigr)}\,du\,=\\
=\,\biggl(\,\int_{h}^{2t}\,+\,\int_{2t}^{+\infty}\,\biggr)K_{a}(\pi u)\ln{\bigl(|t-u|+3\bigr)}\,du\,=\,j_{8}\,+\,j_{9},\\
j_{8}\,\leqslant\,\ln{(t+3)}\int_{h}^{+\infty}K_{a}(\pi u)\,du\,<\,\frac{\ln{(t+3)}}{\pi\ln{T}}\,\frac{1}{\ln\ln{T}},\\
j_{9}\,<\,\int_{2t}^{+\infty}K_{a}(\pi
u)\ln{(u+3)}\,du\,<\,\frac{2}{\pi}\ln{(\pi t)}e^{-\pi
t}\exp{\bigl(-e^{\pi t}\bigr)},
\end{multline*}
and hence $j_{5} < (3\ln\ln{T})^{-1}$.

Going back to (\ref{lab_18}), we obtain
\begin{multline}\label{lab_19}
\int_{-h}^{h}K_{a}(\pi u)\ln{\bigl|\zeta\bigl(0.5+i(t+u)\bigr)\bigr|}\,du\,>\,\frac{2}{a}\,-\\
-\,\biggl(\frac{2}{a}\,e^{-\pi
a/2}\,+\,(7\varepsilon)^{2}\widehat{K}_{a}(0)e^{\pi
a/\sqrt{2}}\,+\,3000\varepsilon^{2}\,+\,\frac{321}{\sqrt{\ln{N}\mathstrut}}\,+\,
\frac{2\widehat{K}_{a}(0)}{3\pi}\,+\,\frac{1}{\ln\ln{T}}\biggr)\,>\\
>\,\frac{2}{a}\,-\,\biggl(\frac{2}{a}\cdot\frac{1}{4}\,+\,(7\varepsilon)^{2}\widehat{K}_{a}(0)e^{\pi a/\sqrt{2}}\,+\,3000\varepsilon^{2}\,
+\,\frac{321}{\sqrt{\ln{N}\mathstrut}}\,+\,\frac{2\widehat{K}_{a}(0)}{3\pi}\biggr).
\end{multline}

In view of (\ref{lab_17}), the expression in the brackets does not exceed
\begin{multline*}
\frac{1}{2a}\,+\,(7\varepsilon)^{2}\sqrt{\frac{2\pi}{a}}\,e^{(\pi/\sqrt{2}-1)a}\,+\,3000\varepsilon^{2}\,\frac{321}{\sqrt{\ln{N}\mathstrut}}\,+\,
\frac{2}{3\pi}e^{-a}\sqrt{\frac{2\pi}{a}}\,<\\
<\,\frac{1}{2a}\biggl(\frac{3}{2}\,+\,2(7\varepsilon)^{2}\sqrt{2\pi
a}\,e^{(\pi/\sqrt{2}-1)a}\,+\,6000a\varepsilon^{2}\,+\,\frac{642a}{\sqrt{\ln{N}\mathstrut}}\biggr).
\end{multline*}
Now we put
\[
\varepsilon\,=\,\frac{e^{-2a/3}}{100\sqrt{a}},\quad
N\,=\,e^{(3852a)^{2}}.
\]
Then the left\,-hand side of the last inequality does not exceed
\[
\frac{1}{2a}\biggl(\frac{3}{2}\,+\,\frac{\sqrt{2\pi}}{100}\,e^{-0.1a}\,+\,\frac{3}{5}\,e^{-4a/3}\,+\,\frac{1}{6}\biggr)\,<\,
\frac{1}{2a}\biggl(\frac{5}{3}\,+\,\frac{1}{6}\,+\,\frac{1}{6}\biggr)\,=\,\frac{1}{a}.
\]
Now (\ref{lab_19}) implies that
\begin{equation}\label{lab_20}
\int_{-h}^{h}K_{a}(\pi
u)\ln{\bigl|\zeta\bigl(0.5+i(t+u)\bigr)\bigr|}\,du\,>\,\frac{2}{a}\,-\,\frac{1}{a}\,=\,\frac{1}{a}.
\end{equation}
Denote by $M$ the maximum of
$\ln{\bigl|\zeta\bigl(0.5+i(t+u)\bigr)\bigr|}$ on the segment
$|u|\leqslant h$. Then (\ref{lab_20}) implies that $M>0$. Hence, the
integral in (\ref{lab_20}) is less than
\[
M\int_{-h}^{h}K_{a}(\pi
u)\,du\,<\,\frac{M}{\pi}\int_{-\infty}^{+\infty}K_{a}(u)\,du\,=\,\frac{M}{\pi}\,\widehat{K}_{a}(0).
\]
Using (\ref{lab_17}) and (\ref{lab_11}), we find
\[
\frac{M}{\pi}\,\widehat{K}_{a}(0)\,>\,\frac{1}{a},\qquad
M\,>\,\frac{\pi}{a}\,\widehat{K}_{a}^{-1}(0)\,>\,e^{a}\sqrt{\frac{\pi}{2a}}\,\geqslant\,\ln{A}.
\]
To end the proof, we note that the distance between $T$ and the point $u$, where the maximum $M$ is attained, does not exceed
\[
c_{0}\,+\,h\,=\,\frac{1}{\pi}\bigl(\ln\ln\ln{T}\,-\,\ln{(a/2)}\bigr)\,+\,c_{0}.
\]
The theorem is proved.

\vspace{0.2cm}

\textsc{Remark.} In \cite{Fyodorov_Keating_2011}, the distribution
of the random variable $\sigma(T)$ with the values
\[
-2\ln{F(t;2\pi)}\,+\,2\ln\ln{\frac{t}{2\pi}}\,-\,\frac{3}{2}\ln\ln\ln{\frac{t}{2\pi}},\quad
t_{0}\leqslant t\leqslant T,
\]
is discussed. In \cite{Harper_2013}, there are some arguments that
reinforce the hypothesis that the inequalities
\[
\frac{\ln{t}}{(\ln\ln{t})^{2+\varepsilon}}\,\leqslant\,F(t;2\pi)\,\leqslant\,\frac{\ln{t}}{(\ln\ln{t})^{0.25-\varepsilon}}
\]
hold for ``almost all'' $t$ from the interval $(T,2T)$, $T\to
+\infty$ and for any $\varepsilon>0$.

\vspace{0.2cm}

\textsc{Theorem 2.} \emph{Suppose that the quantity
\[
S_{0}\,=\,\frac{1}{\pi}\sum\limits_{n\,=\,
p^{2k+1}}(-1)^{k}\,\frac{\Lambda_{1}(n)}{\sqrt{n}}\,\widehat{K}_{a}\biggl(\frac{\ln{n}}{\pi}\biggr)\,=\,\frac{1}{\pi}\Im
\sum\limits_{n
=2}^{+\infty}i^{\,\Omega(n)}\,\frac{\Lambda_{1}(n)}{\sqrt{n}}\,\widehat{K}_{a}\biggl(\frac{\ln{n}}{\pi}\biggr)
\]
is positive for some $a\geqslant 1$. Then for any fixed
$\varepsilon>0$ satisfying the condition
$0<\varepsilon<\min{\bigl(0.5,S_{0}\bigr)}$ there exist the
constants $c_{0}$ and $T_{0}$ depending on $a$ and $\varepsilon$
only and such that the inequalities
\[
\max_{|t-T|\leqslant H}\bigl(\pm
S(t)\bigr)\,>\,\frac{S_{0}-\varepsilon}{\pi \widehat{K}_{a}(0)}
\]
hold for any $T\geqslant T_{0}$ and $H =
(1/\pi)\ln\ln\ln{T}+c_{0}$.}

\vspace{0.2cm}

\textsc{Proof.} Extracting the real parts in (\ref{lab_04}), we obtain:
\begin{equation}\label{lab_21}
\pi\int_{-\infty}^{+\infty}K_{a}(\pi
u)S(t+u)\,du\,=\,-\,\frac{1}{\pi}\sum\limits_{n =
2}^{+\infty}\frac{\Lambda_{1}(n)}{\sqrt{n}}\widehat{K}_{a}\biggl(\frac{\ln{n}}{\pi}\biggr)\sin{(t\ln{n})}\,+\,\frac{\pi\theta_{1}}{\cosh{\pi
t}}.
\end{equation}
Let $\varepsilon_{1}, N$, be the numbers depending on $a$,
$\varepsilon$ and such that $0<\varepsilon_{1}<0.5$, $N\geqslant
e^{\pi a\sqrt{2}}$, whose explicit values will be chosen later. By
Lemma 6, there exists a constant $c = c(a,\varepsilon)$ such that
any interval of length $c$ contains a point $\tau$ such that the
inequality
\begin{equation}\label{lab_22}
\biggl\|\frac{\tau}{2\pi}\ln{p}\,+\,\frac{1}{4}\biggr\|\,<\,\varepsilon_{1}
\end{equation}
holds for any prime $p\leqslant N$. Let us take the parameter $t$ in
(\ref{lab_21}) to be equal to such value $\tau$ from the interval
$(T,T+c)$.

Given prime $p\leqslant N$, we define an integer $n_{p}$ and real
$\varepsilon_{p}$ satisfying the condition
$|\varepsilon_{p}|<\varepsilon_{1}$ such that
\[
\frac{t}{2\pi}\ln{p}\,=\,n_{p}\,+\,\varepsilon_{p}\,-\,\frac{1}{4}.
\]
Then we have
\[
\sin{(t\ln{n})}\,=\,-\sin\left(\frac{\pi k}{2}\right)\cos{(2\pi
k\varepsilon_{p})}\,+\,\cos{\left(\frac{\pi k}{2}\right)}\sin{(2\pi
k\varepsilon_{p})}
\]
for any $k\geqslant 1$ and $n = p^{k}$, $p\leqslant N$. If $k$ is
even then $|\sin{(t\ln{n})}|=|\sin{(2\pi
k\varepsilon_{p}})|\leqslant 2\pi k\varepsilon_{1}$; otherwise, we
have
\[
\sin{(t\ln{n})}\,=\,(-1)^{(k+1)/2}\cos{(2\pi
k\varepsilon_{p})}\,=\,(-1)^{(k+1)/2}\,-\,2\theta_{2}(\pi
k\varepsilon_{1})^{2}.
\]
Let $S$ be the sum in the right\,-hand side of (\ref{lab_21}).
Denote by $S_{1}, S_{2}$ and $S_{3}$ the contributions to this sum
arising from the terms corresponding to the following conditions:
$n= p^{k}$, $p\leqslant N$, $k$ is odd; $n = p^{k}$, $p\leqslant N$,
$k$ is even; $n = p^{k}$, $p> N$, respectively. Then we have
\begin{multline}\label{lab_23}
S_{1}\,=\,-\,\frac{1}{\pi}\sum\limits_{\substack{n\,=\,p^{2k+1} \\
p\,\leqslant\,N,\;k\,\geqslant\,0}}
\frac{\Lambda_{1}(n)}{\sqrt{n}}\,\widehat{K}_{a}\biggl(\frac{\ln{n}}{\pi}\biggr)\Bigl((-1)^{k+1}\,-
\,2\theta_{2}\bigl(\pi(2k+1)\varepsilon_{1}\bigr)^{2}\Bigr)\,=\\
=\,\frac{1}{\pi}\sum\limits_{\substack{n\,=\,p^{2k+1} \\
p\,\leqslant\,N,\;k\,\geqslant\,0}}
(-1)^{k}\,\frac{\Lambda_{1}(n)}{\sqrt{n}}\,\widehat{K}_{a}\biggl(\frac{\ln{n}}{\pi}\biggr)\,+\\
+\,2\theta_{3}\pi\varepsilon_{1}^{2}\sum\limits_{\substack{n\,=\,p^{2k+1}
\\ p\,\leqslant\,N,\;k\,\geqslant\,0}}(2k+1)^{2}\,
\frac{\Lambda_{1}(n)}{\sqrt{n}}\,\biggl|\widehat{K}_{a}\biggl(\frac{\ln{n}}{\pi}\biggr)\biggr|.
\end{multline}
Obviously, the last sum in (\ref{lab_23}) is less than
\begin{multline*}
2\pi\varepsilon_{1}^{2}\widehat{K}_{a}(0)\sum\limits_{p\leqslant
N}\sum\limits_{k =
0}^{+\infty}\frac{2k+1}{p^{k\mathstrut}\sqrt{p}}\,=\,
2\pi\varepsilon_{1}^{2}\widehat{K}_{a}(0)\sum\limits_{p\leqslant N}\frac{1}{\sqrt{p}}\biggl(1\,+\,\frac{1}{p}\biggr)\biggl(1\,-\,\frac{1}{p}\biggr)^{\!-2}\,\leqslant\\
\leqslant\,12\pi\varepsilon_{1}^{2}\widehat{K}_{a}(0)\sum\limits_{p\leqslant
N}\frac{1}{\sqrt{p}}\,<\,\frac{105\varepsilon_{1}^{2}\widehat{K}_{a}(0)\sqrt{N}}{\ln{N}}.
\end{multline*}
Further, we replace the interval $p\leqslant N$ in the first sum in
right\,-hand side of (\ref{lab_23}) by infinite one. The arising
error does not exceed in modulus
\begin{equation}\label{lab_24}
\frac{1}{\pi}\sum\limits_{n>N}\frac{\Lambda_{1}(n)}{\sqrt{n}}\,\biggl|\widehat{K}_{a}\biggl(\frac{\ln{n}}{\pi}\biggr)\biggr|\,<
\,\frac{81}{\sqrt{\ln{N}\mathstrut}}.
\end{equation}
Hence, the difference between $S_{1}$ and $S_{0}$ is less than
\[
\frac{105\varepsilon_{1}^{2}\widehat{K}_{a}(0)\sqrt{N}}{\ln{N}}\,+\,\frac{81}{\sqrt{\ln{N}\mathstrut}}.
\]
Further,
\begin{multline*}
|S_{2}|\,\leqslant\,\frac{1}{\pi}\sum\limits_{\substack{n\,=\,p^{2k} \\
p\leqslant N,\;k\geqslant
1}}\frac{\Lambda_{1}(n)}{\sqrt{n}}\,\bigl|\widehat{K}_{a}\bigl(0\bigr)\bigr|\cdot
4\pi k\varepsilon_{1}\,\leqslant\,
2\varepsilon_{1}\widehat{K}_{a}(0)\sum\limits_{p\leqslant N}\sum\limits_{k = 1}^{+\infty}p^{-k}\,=\\
=\,\,2\varepsilon_{1}\widehat{K}_{a}(0)\sum\limits_{p\leqslant
N}\frac{1}{p-1}\,<\,3\varepsilon_{1}\widehat{K}_{a}(0)\ln\ln{N}.
\end{multline*}
Obviously, the modulus of $S_{3}$ does not exceed the right\,-hand side of (\ref{lab_24}).

Therefore, $S$ and $S_{0}$ differ by at most
\[
\frac{105\varepsilon_{1}^{2}\widehat{K}_{a}(0)\sqrt{N}}{\ln{N}}\,+\,\frac{162}{\sqrt{\ln{N}\mathstrut}}\,+\,3\varepsilon_{1}\widehat{K}_{a}(0)\ln\ln{N}.
\]
We put $h = (1/\pi)\bigl(\ln\ln\ln{T}-\ln{(a/2)}\bigr)$ and split
the improper integral in (\ref{lab_21}) to the integrals $j_{1},
j_{2}$ and $j_{3}$ corresponding to the intervals $|u|\leqslant h$,
$u>h$ and $u<-h$ respectively. If $|v|\geqslant 280$, the classical
Backlund's estimate \cite{Backlund_1916} implies that
\begin{multline}\label{lab_25}
|S(v)|\,<\,0.1361\ln{|v|}\,+\,0.4422\ln\ln{|v|}\,+\,4.3451\,\leqslant\\
\leqslant
\biggl(0.1361\,+\,0.4422\,\frac{\ln\ln{280}}{\ln{280}}\,+\,\frac{4.3451}{\ln{280}}\biggr)\ln{|v|}\,<\,1.05\ln{|v|}.
\end{multline}
Otherwise, we have the inequality $|S(v)|\leqslant 1$ (see
\cite[Tab. 1]{Lehman_1970}). From these estimates, it follows that
$|j_{2}|+|j_{3}|<2(\ln\ln{T})^{-1}$. Hence,
\begin{multline}\label{lab_26}
\pi\int_{-h}^{h}K_{a}(\pi u)S(t+u)\,du\,>\\
>\,S_{0}\,-\,\biggl(\frac{105\varepsilon_{1}^{2}\widehat{K}_{a}(0)\sqrt{N}}{\ln{N}}\,+\,
\frac{162}{\sqrt{\ln{N}\mathstrut}}\,+\,3\widehat{K}_{a}(0)\varepsilon_{1}\ln\ln{N}\,+\,\frac{3}{\ln\ln{T}}\biggr).
\end{multline}
The expression in the brackets is less than
\begin{multline}\label{lab_27}
\frac{105\varepsilon_{1}^{2}\sqrt{N}}{\ln{N}}\,e^{-a}\sqrt{\frac{2\pi}{a}}\,+\,\frac{162}{\sqrt{\ln{N}\mathstrut}}\,+\,
3e^{-a}\sqrt{\frac{2\pi}{a}}\varepsilon_{1}\ln\ln{N}
\,+\,\frac{3}{\ln\ln{T}}\,<\\
<\,\frac{(10\varepsilon_{1})^{2}\sqrt{N}}{\ln{N}}\,+\,\frac{162}{\sqrt{\ln{N}\mathstrut}}\,+\,3\varepsilon_{1}\ln\ln{N}.
\end{multline}
Now we take
\[
\varepsilon_{1}\,=\,\exp{\biggl(-\,\biggl(\frac{162}{\varepsilon}\biggr)^{\!\!2}\biggr)},\qquad
N\,=\,\exp{\biggl(\biggl(\frac{324}{\varepsilon}\biggr)^{\!\!2}\biggr)}.
\]
Then the right\,-hand side of (\ref{lab_27}) is bounded from above by
\[
\biggl(\frac{\varepsilon}{30}\biggr)^{\!
2}\,+\,6\exp{\biggl(-\,\biggl(\frac{162}{\varepsilon}\biggr)^{\!
2}\biggr)}\ln{\biggl(\frac{324}{\varepsilon}\biggr)}\,+\,\frac{\varepsilon}{2}\,<\,\varepsilon.
\]
Since $0<\varepsilon<S_{0}$, the right\,-hand side of (\ref{lab_26})
is positive. Denoting $M_{1} = \displaystyle\max_{|u|\leqslant
h}S(t+u)$, we have therefore:
\[
M_{1}\,>\,0,\qquad S_{0}-\varepsilon\,<\,\pi
M_{1}\int_{-h}^{h}K_{a}(\pi u)\,du\,<\,M_{1}\widehat{K}_{a}(0).
\]
Thus, $M_{1}>(S_{0}-\varepsilon)\widehat{K}_{a}^{-1}(0)$. Since the
distance between $T$ and the point $t+u$, where the maximum is
attained, is less than $H = h+c =
(1/\pi)\bigl(\ln\ln\ln{T}-\ln{(a/2)}\bigr)+c$, the first statement
of theorem is proved. The proof of second one is similar. The only
difference is that $t$ is chosen now in $(T,T+c)$ to satisfy the
inequalities
\[
\biggl\|\frac{t}{2\pi}\ln{p}\,-\,\frac{1}{4}\biggr\|\,<\,\varepsilon_{1}
\]
for all primes $p\leqslant N$. The theorem is proved.

\vspace{0.2cm}

The very slow convergence of the series $S_{0}$ and the absence of the analogue of the identity (\ref{lab_13})
make the verification of the condition $S_{0}>0$ very difficult. However, a small modification of the
above proof allows one to obtain a series of numerical results.

\vspace{0.2cm}

\textsc{Theorem 3.} \emph{Let $a, b, \tau$ be any positive numbers
satisfying the conditions $0<b<\pi/2$, $b\tau>0.5$, $\gamma =
b\tau+0.5$, $N\geqslant 2$ be an integer, and let}
\[
S_{N}(u)\,=\,\sum\limits_{p\leqslant
N}\arctan{\biggl(\frac{2\sqrt{p}}{p-1}\cos{(u\tau\ln{p})}\biggr)}.
\]
\emph{Further, let}
\begin{align*}
& \kappa\,=\,\kappa(a,b)\,=\,2\int_{0}^{+\infty}e^{-a\cos{(b)}\cosh{(u)}}\,du,\\
&
\zeta_{N}(\gamma)\,=\,\prod\limits_{p>N}\bigl(1\,-\,p^{-\,\gamma}\bigr)^{-1}\,=\,\zeta(\gamma)\prod\limits_{p\leqslant
N}\bigl(1\,-\,p^{-\,\gamma}\bigr),
\end{align*}
\emph{and let}
\[
I\,=\,\frac{1}{\pi}\biggl(\int_{0}^{+\infty}K_{a}(u)S_{N}(u)\,du\,-\,\kappa\ln{\zeta_{N}(\gamma)}\biggr)\,>\,0.
\]
\emph{Then, for any fixed $\varepsilon$,
$0<\varepsilon<\varepsilon_{0}(a,b,\tau)$, there exists a constant
$c_{0} = c_{0}(\varepsilon; a,b,\tau)$ such that the inequalities}
\[
\max_{|T-t|\leqslant H}\bigl(\,\pm
S(t)\,\bigr)\,>\,\frac{I-\varepsilon}{\widehat{K}_{a}(0)}
\]
\emph{hold for any $T\geqslant T_{0}(\varepsilon; a,b,\tau)$ and $H
= \tau\ln\ln\ln{T}+c_{0}$}.

\vspace{0.2cm}

\textsc{Proof.} Setting $f(u) = (1/\tau)K_{a}(u/\tau)$ in Lemma 4
and extracting imaginary parts, we get
\begin{equation}\label{lab_28}
\frac{1}{\tau}\int_{-\infty}^{+\infty}K_{a}\biggl(\frac{u}{\tau}\biggr)S(t+u)\,du\,=\,C\,+\,\frac{\pi\theta_{1}}{a\cosh{(t/\tau)}},
\end{equation}
where
\[
C\,=\,-\,\frac{1}{\pi}\sum\limits_{n =
2}^{+\infty}\frac{\Lambda_{1}(n)}{\sqrt{n}}\,\widehat{K}_{a}(\tau\ln{n})\sin{(\tau\ln{n})}.
\]
We denote
\[
c\,=\,\kappa\biggl(\frac{4}{3}+7\ln{\zeta(\gamma)}\biggr),\quad
\varepsilon_{0}\,=\,\min{\biggl(0.5,\frac{I}{c},\frac{\varepsilon}{c}\biggr)}
\]
and take an arbitrary fixed numbers $\varepsilon_{1}$ and $N$
satisfying the conditions $0<\varepsilon_{1}<\varepsilon_{0}$,
$N>2$. By Corollary of Lemma 6, there exists a constant $c_{0}$
depending only on $\varepsilon_{1}$, $N$ and such that any interval
of length $c_{0}$ contains a point $\tau$ such that the inequality
(\ref{lab_22}) holds for any prime $p\leqslant N$. Suppose that $t$
is such value from the interval $(T,T+c_{0})$. Similarly to the
proof of Theorem 2, we split the sum $C$ to the sums $C_{1}, C_{2}$
and $C_{3}$. Thus we get $C_{1} = C_{0} + \theta_{2}C_{4}$, where
\begin{multline*}
C_{0}\,=\,\frac{1}{\pi}\sum\limits_{\substack{n = p^{2k+1\mathstrut},\;k\geqslant 0 \\ p\leqslant N}}(-1)^{k}\,\frac{\Lambda_{1}(n)}{\sqrt{n}}\,\widehat{K}_{a}(\tau\ln{n}),\\
C_{4}\,=\,2\pi\varepsilon_{1}^{2}\sum\limits_{\substack{n =
p^{2k+1\mathstrut},\;k\geqslant 0 \\ p\leqslant
N}}(2k+1)^{2}\,\frac{\Lambda_{1}(n)}{\sqrt{n}}\,\bigl|\widehat{K}_{a}(\tau\ln{n})\bigr|.
\end{multline*}
Moreover,
\[
|C_{2}|\,\leqslant\,4\varepsilon_{1}\sum\limits_{\substack{n =
p^{2k\mathstrut},\;k\geqslant 1 \\ p\leqslant
N}}k\,\frac{\Lambda_{1}(n)}{\sqrt{n}}\,\bigl|\widehat{K}_{a}(\tau\ln{n})\bigr|,\quad
|C_{3}|\,\leqslant\,\frac{1}{\pi}\sum\limits_{n=p^{k},\;
p>N}\frac{\Lambda_{1}(n)}{\sqrt{n}}\,\bigl|\widehat{K}_{a}(\tau\ln{n})\bigr|.
\]
The application of Lemma 2 yields:
\begin{multline*}
|C_{2}|\,\leqslant\,4\varepsilon_{1}\sum\limits_{\substack{n = p^{2k\mathstrut},\;k\geqslant 1 \\ p\leqslant N}}\frac{k}{2k\sqrt{n}}\,\kappa e^{-\,b\tau\ln{n}}\,=\,2\kappa\varepsilon_{1}\sum\limits_{p\leqslant N}\sum\limits_{k = 1}^{+\infty}p^{-\,2k\gamma}\,=\\
=\,2\kappa\varepsilon_{1}\sum\limits_{p\leqslant
N}p^{-2\gamma}\bigl(1\,-\,p^{-\,2\gamma}\bigr)^{-1}\,\leqslant\,\frac{2\kappa\varepsilon_{1}}{1-2^{-2\gamma}}\,\ln{\zeta(2\gamma)}\,
<\,\frac{2\kappa\varepsilon_{1}}{1-2^{-2}}\,\ln{\zeta(2)}\,<\,\frac{4}{3}\,\kappa\varepsilon_{1},\\
|C_{3}|\,\leqslant\,\frac{\kappa}{\pi}\sum\limits_{\substack{n =
p^{k\mathstrut},\;k\geqslant 1 \\ p >
N}}\frac{\Lambda_{1}(n)}{n^{\gamma}}\,=\,\frac{\kappa}{\pi}\sum\limits_{p>N}\ln{\bigl(1\,-\,p^{-\,\gamma}\bigr)^{-1}}\,=\,
\frac{\kappa}{\pi}\ln{\biggl(\zeta(\gamma)\prod\limits_{p\leqslant N}\bigl(1\,-\,p^{-\,\gamma}\bigr)\biggr)}\,=\\
=\,\frac{\kappa}{\pi}\ln{\zeta_{N}(\gamma)},
\end{multline*}
and finally
\begin{multline*}
|C_{4}|\,\leqslant\,2\pi\kappa\varepsilon_{1}^{2}\sum\limits_{\substack{n
= p^{2k+1\mathstrut},\;k\geqslant 0 \\ p\leqslant
N}}\frac{2k+1}{n^{\gamma}}\,=\,
2\pi\kappa\varepsilon_{1}^{2}\sum\limits_{p\leqslant N}\sum\limits_{k = 0}^{+\infty}\frac{2k+1}{p^{(2k+1)\gamma}}\,=\\
=\,2\pi\kappa\varepsilon_{1}^{2}\sum\limits_{p\leqslant
N}\frac{1}{p^{\,\gamma}}\,\frac{1+p^{-2\gamma}}{(1\,-\,p^{-\,2\gamma})^{2}}\,
\leqslant\,2\pi\kappa\varepsilon_{1}^{2}\,\frac{1+2^{-2\gamma}}{(1\,-\,2^{-2\gamma})^{2}}\sum\limits_{p\leqslant
N}\frac{1}{p^{\,\gamma}}\,<\,\frac{40\pi}{9}\kappa\varepsilon_{1}^{2}\ln{\zeta(\gamma)}.
\end{multline*}
Transforming the sum $C_{0}$, we obtain
\begin{multline}\label{lab_29}
C_{0}\,=\,\frac{1}{\pi}\sum\limits_{\substack{n =
p^{2k+1\mathstrut},\;k\geqslant 0 \\ p\leqslant
N}}(-1)^{k}\,\frac{\Lambda_{1}(n)}{\sqrt{n}}
\int_{-\infty}^{+\infty}K_{a}(u)e^{-\,iu\tau\ln{n}}\,du\,=\\
=\,\frac{1}{\pi}\int_{-\infty}^{+\infty}K_{a}(u)
\biggl(\;\;\sum\limits_{\substack{n = p^{2k+1\mathstrut},\;k\geqslant 0 \\ p\leqslant N}}(-1)^{k}\,\frac{\Lambda_{1}(n)}{\sqrt{n}}\,n^{-iu\tau}\biggr)du\,=\\
=\,\frac{1}{\pi}\int_{-\infty}^{+\infty}K_{a}(u)\sum\limits_{p\leqslant
N}\biggl(\,\sum\limits_{k = 0}^{+\infty}
\frac{(-1)^{k}}{2k+1}\biggl(\frac{p^{-\,iu\tau}}{\sqrt{p}}\biggr)^{\!\! 2k+1}\,\biggr)\,du\,=\\
=\,\frac{1}{2\pi i}\int_{0}^{+\infty}K_{a}(u)\sum\limits_{p\leqslant
N}\biggl\{\ln{\biggl(1+\frac{ip^{-iu\tau}}{\sqrt{p}}\biggr)}\,-
\,\ln{\biggl(1-\frac{ip^{-iu\tau}}{\sqrt{p}}\biggr)}\,+\\
+\,\ln{\biggl(1+\frac{ip^{iu\tau}}{\sqrt{p}}\biggr)}\,-\,\ln{\biggl(1-\frac{ip^{iu\tau}}{\sqrt{p}}\biggr)}\biggr\}\,du.
\end{multline}
For fixed $p\leqslant N$, we denote
\[
z_{1}\,=\,1\,+\,\frac{ip^{iu\tau}}{\sqrt{p}}\,=\,|z_{1}|e^{i\varphi_{1}},\qquad
z_{2}\,=\,1\,-\,\frac{ip^{iu\tau}}{\sqrt{p}}\,=\,|z_{2}|e^{i\varphi_{2}},
\]
where $-\pi <\varphi_{1},\varphi_{2}\leqslant \pi$. Then the
summands in (\ref{lab_29}) take the form
\[
\ln{\overline{z}_{2}}\,-\,\ln{\overline{z}_{1}}\,+\,\ln{z_{1}}\,-\,\ln{z_{2}}\,=\,
\ln{\frac{\overline{z}_{2}}{z_{2}}}\,-\,\ln{\frac{\overline{z}_{1}}{z_{1}}}
\,=\,2i(\varphi_{1}\,-\,\varphi_{2}).
\]
Writing $\alpha_{p} = u\tau\ln{p}$ and noting that
\[
z_{1}\,=\,1\,-\,\frac{\sin{\alpha_{p}}}{\sqrt{p}}\,+\,\frac{i\cos{\alpha_{p}}}{\sqrt{p}},
\]
we find
\[
\tan{\varphi_{1}}\,=\,\tan{(\arg{z_{1}})}\,=\,\frac{(\cos{\alpha_{p}})/\sqrt{p}}{1-(\sin{\alpha_{p}})/\sqrt{p}}\,=
\,\frac{\cos{\alpha_{p}}}{\sqrt{p}-\sin{\alpha_{p}}}.
\]
Similarly,
\[
\tan{\varphi_{2}}\,=\,\tan{(\arg{z_{2}})}\,=\,-\,\frac{\cos{\alpha_{p}}}{\sqrt{p}+\sin{\alpha_{p}}}.
\]
Hence
\[
\tan{(\varphi_{1} -
\varphi_{2})}\,=\,\frac{\tan{\varphi_{1}}\,-\,\tan{\varphi_{2}}}{1\,+\,\tan{\varphi_{1}}\tan{\varphi_{2}}}\,=\,\frac{2\sqrt{p}}{p-1}\,\cos{\alpha_{p}},
\]
and therefore
\begin{multline*}
\varphi_{1}\,-\,\varphi_{2}\,=\,\arctan{\biggl(\frac{2\sqrt{p}}{p-1}\,\cos{\alpha_{p}}\biggr)},\\
C_{0}\,=\,\frac{1}{\pi}\int_{0}^{+\infty}K_{a}(u)\sum\limits_{p\leqslant
N}\arctan{\biggl(\frac{2\sqrt{p}}{p-1}\,\cos{\alpha_{p}}\biggr)}\,du\,=\,
\frac{1}{\pi}\int_{0}^{+\infty}K_{a}(u)S_{N}(u)\,du.
\end{multline*}
Summing the above bounds, we conclude that the difference between
$C_{0}$ and the right\,-hand side of (\ref{lab_28}) does not exceed
in modulus
\begin{multline*}
\frac{\kappa}{\pi}\ln{\zeta_{N}(\gamma)}\,+\,\frac{4}{3}\,\kappa\varepsilon_{1}\,+\,\frac{40\pi}{9}\kappa\varepsilon_{1}^{2}\ln{\zeta(\gamma)}\,+\,
\frac{\pi}{a\cosh{(t/\tau)}}\,<\\
<\,\frac{\kappa}{\pi}\ln{\zeta_{N}(\gamma)}\,+\,\kappa\varepsilon_{1}\biggl(\frac{4}{3}\,+\,7\ln\zeta(\gamma)\biggr)
-\frac{3}{\ln\ln{T}}\,=\,\frac{\kappa}{\pi}\ln{\zeta_{N}(\gamma)}\,+\,c\varepsilon_{1}\,-\,\frac{3}{\ln\ln{T}}.
\end{multline*}
Let $h = \tau\bigl(\ln\ln\ln{T}-\ln{(a/2)}\bigr)$. Splitting the
integral in (\ref{lab_28}) to the sum
\[
j_{1}\,+\,j_{2}\,+\,j_{3}\,=\,\frac{1}{\tau}\biggl(\;\int_{-h}^{h}\,+\,\int_{h}^{+\infty}\,+\,
\int_{-\infty}^{-h}\;\biggr)K_{a}\biggl(\frac{u}{\tau}\biggr)S(t+u)\,du
\]
and using the same bounds for $S(u)$ as in the proof of Theorem 2,
we find: $|j_{2}|+|j_{3}|<3(\ln\ln{T})^{-1}$. Hence,
\[
j_{1}\,=\,\frac{1}{\tau}\int_{-h}^{h}K_{a}\biggl(\frac{u}{\tau}\biggr)S(t+u)\,du\,>\,C_{0}\,-\,\frac{\kappa}{\pi}\ln{\zeta_{N}(\gamma)}\,-\,
c\varepsilon_{1}\,=\,I\,-\,c\varepsilon_{1}.
\]
Since $0<\varepsilon_{1}<I/c$, the right\,-hand side of the last
inequality is strictly positive, and so is the quantity $M_{1} =
\displaystyle\max_{|u|\leqslant h}S(t+u)$. Obviously, we have
$j_{1}<M_{1}\widehat{K}_{a}(0)$, and therefore
$M_{1}>(I-\varepsilon)/\widehat{K}_{a}(0)$. The lower bound of
$M_{2} = \displaystyle\max_{|u|\leqslant h}\bigl(-\,S(t+u)\bigr)$ is
established by similar arguments. The theorem is proved.

\vspace{0.5cm}

The condition $I>0$ can be checked without significant difficulties.
Let
\[
\mu\,=\,\frac{I}{\widehat{K}_{a}(0)}\,=\,
\frac{1}{\pi\widehat{K}_{a}(0)}\biggl(\;\int_{0}^{+\infty}K_{a}(u)S_{N}(u)\,du\,-\,\kappa\ln{\zeta_{N}(\gamma)}\biggr).
\]
Taking $a = 3$, $b = 7/5$, $\tau = 2/5$ and choosing $N = p_{n}$
from the table below, we find that
\begin{align*}
& n && \mu \\
& 16\,500 && 1.005\,075\,13\ldots \\
& 78\,000 && 2.006\,322\,98\ldots \\
& 2\,500\,000 && 3.001\,263\,70\ldots \;.
\end{align*}

\textsc{Corollary.} \emph{If the Riemann hypothesis is true, than
there exist the constants $c_{0}$ and $T_{0}$ such that the
inequalities}
\[
\max_{|t-T|\leqslant H}\bigl(\pm S(t)\bigr)\,>\,3\,+\,10^{-3}
\]
\emph{hold for any} $T\geqslant T_{0}$ and $H =
0.4\ln\ln\ln{T}+c_{0}$.

\vspace{0.2cm}

\textsc{Theorem 4.} \emph{Suppose that the Riemann hypothesis is
true. Then for an arbitrary large fixed $A\geqslant 1$, there exist
constants $T_{0}, c_{0}$ and $h$ depending only on $A$ and such that
the inequality
\[
\min_{|t-T|\leqslant H}\bigl(S(t+h)\,-\,S(t-h)\bigr)\,<\,-A
\]
holds with $T\geqslant T_{0}$ and $H = (1/\pi)\ln\ln\ln{T}+c_{0}$.}

\vspace{0.2cm}

\textsc{Proof.} Let $a>1$ and $0<h<1$ be fixed numbers. Replacing $t$ in (\ref{lab_04}) by $t+h$ and $t-h$ and
subtracting the corresponding relations, we obtain
\begin{multline}\label{lab_30}
\int_{-\infty}^{+\infty}K_{a}(\pi u)\bigl(\ln{\zeta(0.5+i(t+h))}\,-\,\ln{\zeta(0.5+i(t-h))}\bigr)\,du\,=\\
=\,\frac{2}{\pi i}\sum\limits_{n = 2}^{+\infty}\frac{\Lambda_{1}(n)}{\sqrt{n}}\,\widehat{K}_{a}\biggl(\frac{\ln{n}}{\pi}\biggr)\sin{(h\ln{n})}n^{-it}\,-\\
-\,2\pi\int_{0}^{0.5}\bigl(K_{a}(\pi(t+h+iv))\,-\,K_{a}(\pi(t-h+iv))\bigr)\,dv.
\end{multline}
Taking imaginary parts in (\ref{lab_30}), we get
\begin{multline}\label{lab_31}
\pi\int_{-\infty}^{+\infty}K_{a}(\pi u)\bigl(S(t+h+u)\,-\,S(t-h+u)\bigr)\,du\,=\\
=\,-\,\frac{2}{\pi}\sum\limits_{n = 2}^{+\infty}\frac{\Lambda_{1}(n)}{\sqrt{n}}\,\widehat{K}_{a}\biggl(\frac{\ln{n}}{\pi}\biggr)\sin{(h\ln{n})}\cos{(t\ln{n})}\,-\\
-\,2\pi\Im\int_{0}^{0.5}\bigl(K_{a}(\pi(t+h+iv))\,-\,K_{a}(\pi(t-h+iv))\bigr)\,dv.
\end{multline}
If $t = 0$ then the integral in the right\,-hand side in (\ref{lab_31}) has the form
\begin{multline*}
-2\pi\Im\int_{0}^{0.5}e^{-a\cosh{(\pi h)}\cos{(\pi v)}}\bigl(\,e^{-ia\sinh{(\pi h)}\sin{(\pi v)}}\,-\,e^{ia\sinh{(\pi h)}\sin{(\pi v)}}\,\bigr)\,dv\,=\\
=\,4\pi\int_{0}^{0.5}e^{-a\cosh{(\pi h)}\cos{(\pi
v)}}\sin{\bigl(a\sinh{(\pi h)}\sin{(\pi v)}\bigr)}\,dv.
\end{multline*}
Hence, we have
\begin{multline}\label{lab_32}
-\,\frac{2}{\pi}\sum\limits_{n = 2}^{+\infty}\frac{\Lambda_{1}(n)}{\sqrt{n}}\,\widehat{K}_{a}\biggl(\frac{\ln{n}}{\pi}\biggr)\sin{(h\ln{n})}\,=\\
=\,-\,4\pi\int_{0}^{0.5}e^{-a\cosh{(\pi h)}\cos{(\pi v)}}\sin{\bigl(a\sinh{(\pi h)}\sin{(\pi v)}\bigr)}\,dv\,+\\
+\,\pi\int_{-\infty}^{+\infty}K_{a}(\pi
u)\bigl(S(u+h)\,-\,S(u-h)\bigr)\,du.
\end{multline}
Let $\varepsilon, N$ be the numbers satisfying the conditions
$0<\varepsilon<0.5$, $N>e^{\pi a\sqrt{2}}$ and depending only on
$a$, whose precise values will be chosen below.

By Lemma 6, given $\varepsilon, N$ satisfying the conditions
$0<\varepsilon<0.5$, $N>e^{\pi a\sqrt{2}}$, there exists a constant
$c$ such that any interval of length $c$ contains a point $\tau$
such that $\bigl\|(\tau/(2\pi))\ln{p}\bigr\|<\varepsilon$ for any
prime $p\leqslant N$. Taking $t$ in (\ref{lab_31}) to be equal to
such value from the interval $(T,T+c)$, estimating the integral in
the right\,-hand side of (\ref{lab_31}) by
$2\pi\bigl(a\cosh{\pi(t-h)}\bigr)^{-1}$ and using the identity
(\ref{lab_32}), we transform the right\,-hand side of (\ref{lab_31})
to the form
\begin{multline}\label{lab_33}
-\,4\pi\int_{0}^{0.5}e^{-a\cosh{(\pi h)}\cos{(\pi v)}}\sin{\bigl(a\sinh{(\pi h)}\sin{(\pi v)}\bigr)}\,dv\,+\\
+\,\pi\int_{-\infty}^{+\infty}K_{a}(\pi u)\bigl(S(u+h)\,-\,S(u-h)\bigr)\,du\,+\\
+\,\frac{4}{\pi}\sum\limits_{n =
2}^{+\infty}\frac{\Lambda_{1}(n)}{\sqrt{n}}\,\widehat{K}_{a}\biggl(\frac{\ln{n}}{\pi}\biggr)\sin{(h\ln{n})}\sin^{2}{\biggl(\frac{t}{2}\ln{n}\biggr)}
\,+\,\frac{2\pi\theta_{1}}{a\cosh{\pi(t-h)}}.
\end{multline}
The sum over $n$ in the right\,-hand side of (\ref{lab_33}) is estimated in the same way as the sum $C$ in Theorem 1 and does not exceed
\[
2\biggl((7\varepsilon)^{2}\widehat{K}_{a}(0)e^{\pi
a/\sqrt{2}}\,+\,3000\varepsilon^{2}\,+\,\frac{321}{\sqrt{\ln{N\mathstrut}}}\biggr)
\]
in modulus. In view of (\ref{lab_25}), the improper integral in (\ref{lab_32}) does not exceed
\begin{multline*}
2\pi\int_{-279}^{279}K_{a}(\pi
u)du\,+\,\pi\biggl(\,\int_{279}^{+\infty}\,+\,\int_{-\infty}^{-279}\,\biggr)K_{a}(\pi
u)\cdot
2.1\ln{(|u|+1)}du\,<\\
<\,2\widehat{K}_{a}(0)\,+\,10^{-100}\widehat{K}_{a}(0)\,<\,2.1\widehat{K}_{a}(0)
\end{multline*}
in absolute value. Hence, changing the signs in (\ref{lab_33}), we get
\begin{multline}\label{lab_34}
\pi\int_{-\infty}^{+\infty}K_{a}(\pi u)\bigl(S(u+h)\,-\,S(u-h)\bigr)\,du\,>\\
>\,4\pi\int_{0}^{0.5}e^{-a\cosh{(\pi h)}\cos{(\pi v)}}\sin{\bigl(a\sinh{(\pi h)}\sin{(\pi v)}\bigr)}\,dv\,-\\
-\,2\biggl((7\varepsilon)^{2}\widehat{K}_{a}(0)e^{\pi
a/\sqrt{2}}\,+\,3000\varepsilon^{2}\,+\,\frac{3210}{\sqrt{\ln{N\mathstrut}}}\,+\,
2.1\widehat{K}_{a}(0)\,+\,\frac{2\pi}{a\cosh{\pi(t-h)}}\biggr).
\end{multline}
Now we take $h = (2\pi a)^{-1}$ and estimate the integral in the right\,-hand side of (\ref{lab_34}) from below. Since
\[
\sin{\bigl(a\sinh{(\pi h)}\sin{(\pi
v)}\bigr)}\,\geqslant\,\sin{\biggl(a\pi h\cdot\frac{2}{\pi}\,\pi
v\biggr)}\,=\,\sin{v}\,\geqslant\,\frac{2}{\pi}\,v, \quad \cosh{\pi
h}< \cosh{\frac{1}{2}}<\frac{8}{7},
\]
the integral under considering is greater than
\begin{multline*}
4\pi\int_{0}^{0.5}e^{-(8a/7)\cos{(\pi v)}}\frac{2}{\pi}v\,dv\,=\,\frac{8}{\pi^{2}}\int_{0}^{\pi/2}e^{-(8a/7)\cos{w}}\,w\,dw\,=\\
=\,\frac{8}{\pi^{2}}\int_{0}^{\pi/2}e^{-(8a/7)\sin{w}}\biggl(\frac{\pi}{2}\,-\,w\biggr)\,dw\,
\geqslant\,\frac{2}{\pi}\int_{0}^{\pi/4}e^{-(8a/7)\sin{w}}\,dw\,\geqslant\\
\geqslant\,\frac{2}{\pi}\int_{0}^{\pi/4}e^{-(8a/7)w}\,dw\,=\,\frac{7}{4\pi
a}\bigl(1\,-\,e^{-2\pi a/7}\bigr)\,>\,\frac{7}{4\pi
a}\bigl(1\,-\,e^{-2\pi/7}\bigr)\,>\,\frac{0.33}{a}.
\end{multline*}
Therefore,
\begin{multline*}
\pi\int_{-\infty}^{+\infty}K_{a}(\pi u)\bigl(S(t+u-h)\,-\,S(t+u+h)\bigr)\,du\,>\,\frac{0.33}{a}\,-\\
-\,\biggl(2(7\varepsilon)^{2}\widehat{K}_{a}(0)e^{\pi
a/\sqrt{2}}\,+\,3000\varepsilon^{2}\,+\,\frac{642}{\sqrt{\ln{N\mathstrut}}}\,+\,2.1\widehat{K}_{a}(0)\,+\,\frac{2\pi}{a\cosh{\pi(t-h)}}\biggr).
\end{multline*}
Let $H_{0} = (1/\pi)\bigl(\ln\ln\ln{T}-\ln{(a/2)}\bigr)$. Then the
sum of integrals over the intervals $(-\infty,-H_{0})$ and
$(H_{0},+\infty)$ in the right\,-hand side is less than
$(\ln\ln{T})^{-1}$ in modulus. Thus we get
\begin{multline}\label{lab_35}
\pi\int_{-H_{0}}^{H_{0}}K_{a}(\pi u)\bigl(S(t+u-h)\,-\,S(t+u+h)\bigr)\,du\,>\,\frac{0.33}{a}\,-\\
-\,\biggl(2(7\varepsilon)^{2}\widehat{K}_{a}(0)e^{\pi
a/\sqrt{2}}\,+\,6000\varepsilon^{2}\,+\,\frac{642}{\sqrt{\ln{N\mathstrut}}}\,+\,2.1\widehat{K}_{a}(0)\,+\,\frac{2}{\ln\ln{T}}\biggr).
\end{multline}
Suppose now that $a>8$ and take $\varepsilon =
e^{-2a/3}/(65\sqrt{a})$, $N = e^{(c_{1}a)^{2}}$, $c_{1}=2^{16}$.
Then
\begin{multline*}
(2(7\varepsilon)^{2}\widehat{K}_{a}(0)e^{\pi
a/\sqrt{2}}\,+\,6000\varepsilon^{2}\,<\,98\varepsilon^{2}e^{-a}\sqrt{\frac{2\pi}{a}}\,e^{\pi
a/\sqrt{2}}\,+\,6000\varepsilon^{2}\,<\\
<\,\frac{98\sqrt{2\pi}}{65^{2}}\,\frac{e^{-0.1a}}{\sqrt{a}}\,\frac{1}{a}\,+\,
\frac{6000}{65^{2}}\,\frac{e^{-4a/3}}{a}\,<\,\frac{10^{-2}}{a},\\
\frac{642}{\sqrt{\ln{N}}}\,=\,\frac{642}{2^{16}a}\,<\,\frac{10^{-2}}{a},\qquad
2.1\widehat{K}_{a}(0)\,+\,\frac{2}{\ln\ln{T}}\,<\,2.1e^{-a}\sqrt{2\pi
a}\,\frac{1}{a}\,<\,\frac{5\cdot 10^{-3}}{a}.
\end{multline*}

Thus, the right\,-hand side of (\ref{lab_35}) is bounded from below by
\[
\frac{0.33}{a}\,-\,\biggl(\frac{2\cdot 10^{-2}}{a}\,+\,\frac{5\cdot
10^{-3}}{a}\biggr)\,>\,\frac{0.3}{a}.
\]
Hence, the value
\[
M_{0}=\displaystyle \max_{|u|\leqslant
H_{0}}\bigl(S(t+u-h)-S(t+u+h)\bigr)
\]
is positive, and the left-\,hand side of (\ref{lab_35}) does not exceed $M_{0}\widehat{K}_{a}(0)$.
Therefore,
\[
M_{0}\,>\,\frac{3\widehat{K}_{a}^{-1}(0)}{10a}\,>\,\frac{3e^{a}}{10\sqrt{2\pi
a}}\,>\,\frac{e^{a}}{10\sqrt{a}}.
\]
Choosing $a>8$ such that
\[
\frac{e^{a}}{10\sqrt{a}}\,>\,A,
\]
we arrive at the assertion of the theorem. The theorem is proved.

In  \cite{Mueller_1983},
\cite{Balasubramanian_1986}, \cite{Korolev_2005},
\cite{Karatsuba_2006} и \cite{Boyarinov_2010}, one can find some other
examples of application the function $K_{a}(z)$ to the theory of $\zeta(s)$.

The key ingredient of the proof of the unboundedness of
$\bigl|\zeta(0.5+it)\bigr|$ on the segment $|t-T|\ll \ln\ln\ln{T}$
is the presence of the term
\[
2\pi\int_{0}^{0.5}e^{-\,a\cos{(\pi v)}}\,dv
\]
in the right\,-hand side of (\ref{lab_13}). It follows from the proof of (\ref{lab_03}) that the pole of $\zeta(s)$ at
the point $s = 1$ is the reason of the appearance of that term.
In view of this, it is interesting to prove the analogue of Theorem 1 for the functions that are ``similar'' to $\zeta(s)$ but have no
pole at the point $s = 1$ (for example, for Dirichlet's $L$\,-function $L(s,\chi_{4})$, where $\chi_{4}$ is non\,-principal character mod $4$).

\vspace{0.5cm}

\textbf{\S 4. The distribution of zeros of zeta\,-function.}

\vspace{0.5cm}

The above theorems allow one to establish some new statements concerning the distribution
of zeros of the Riemann zeta function. Here we also suppose that the Riemann hypothesis is true.

Let $N(t)$ be the number of zeros of $\zeta(s)$ whose ordinate is positive and does not exceed $t$. Then it is known that
\[
N(t)\,=\,\frac{1}{\pi}\,\vartheta(t)\,+\,1\,+\,S(t)\,=\,\frac{t}{2\pi}\ln\frac{t}{2\pi}\,-\,\frac{t}{2\pi}\,+\,\frac{7}{8}\,+\,S(t)\,+\,O\bigl(t^{-1}\bigl),
\]
where $\vartheta(t)$ denotes the increment of a continuous branch of
the argument of the function $\pi^{-s/2}\Gamma(s/2)$ along the line
segment joining the points $s = 0.5$ and $s = 0.5+it$. Then the
Gram's point $t_{n}$ ($n\geqslant 0$) is defined as a unique
solution of the equation $\vartheta(t_{n}) = (n-1)\pi$ with the
condition $\vartheta'(t_{n})>0$. It is easy to check that the number
of zeros of $\zeta\bigl(0.5+it\bigr)$ lying in the Gram's interval
$G_{n} = (t_{n-1}, t_{n}]$ is equal to
\begin{equation}\label{lab_36}
N(t_{n}+0)\,-\,N(t_{n-1}+0)\,=\,1\,+\,\Delta(n)\,-\,\Delta(n-1),
\end{equation}
where $\Delta(n) = S(t_{n}+0)$. Since the segment $[0,T]$ contains
\[
\frac{1}{\pi}\vartheta(T)\,+\,O(1)\,=\,N(T)\,+\,O(\ln{T})
\]
Gram's intervals $G_{n}$, there is precisely one zero of
$\zeta\bigl(0.5+it\bigr)$ per one Gram's interval $G_{n}$ ``in the
mean''. That is the reason why the difference $\Delta(n) -
\Delta(n-1)$ in (\ref{lab_36}) is the deviation of number of zeros
of $\zeta\bigl(0.5+it\bigr)$ in the interval $G_{n}$ from its mean
value, that is, $1$.

In 1946, A.~Selberg \cite{Selberg_1946b} proved that the interval $G_{n}$ contains no zeros of $\zeta\bigl(0.5+it\bigr)$ for positive proportion of $n$,
and contains at least two zeros for positive proportion of $n$ at the same time. These facts show the evident irregularity in the distribution of zeta zeros.

However, nothing is known about the distribution of Gram's intervals $G_{n}$ which are ``free'' of zeros of $\zeta\bigl(0.5+it\bigr)$.
The below theorem establishes an upper bound for the length $h = h(t)$ of the interval $(t,t+h)$ which certainly contains an ``empty'' Gram's interval
$G_{n}$.

\vspace{0.2cm}

\textsc{Theorem 5.} \emph{Suppose that the Riemann hypothesis is
true and let $\varepsilon$ be any fixed positive constant. Then
there exist constants $T_{0}=T_{0}(\varepsilon)$ and $c_{0} =
c_{0}(\varepsilon)$ such that any segment $[T-H,T+H]$, where
$T\geqslant T_{0}$ and $H = (1/\pi)\ln\ln\ln{T}+c_{0}$, contains at
least
$N=\bigl[0.1\sqrt{\varepsilon}\exp{\bigl((\pi\varepsilon)^{-1}\bigr)}\bigr]$
Gram's intervals $G_{n} = (t_{n-1},t_{n}]$ that do not contain zeros
of $\zeta\bigl(0.5+it\bigr)$. Moreover, there exist at least $N$
intervals among the above ``empty'' Gram's intervals that lie in the
same segment of length $\varepsilon$.}

\vspace{0.2cm}

\textsc{Proof.} Let $a = (\pi\varepsilon)^{-1}$, $h = (2\pi a)^{-1}
= 0.5\varepsilon$ and suppose $\varepsilon$ to be so small that $M =
e^{a}/(10\sqrt{a})\geqslant 5$. By Theorem 4, there exist constants
$T_{0} = T_{0}(\varepsilon)$ and $c_{1} = c_{1}(\varepsilon)$ such
that the inequality
\begin{equation*}
\min_{|t-T|\leqslant
H}\bigl(\,S(t+h)\,-\,S(t-h)\,\bigr)\,\leqslant\,-M
\end{equation*}
holds for any $T\geqslant T_{0}$ with $H =
(1/\pi)\ln\ln\ln{T}+c_{1}$.

Let $k$ be sufficiently large and suppose that $t_{k-1}\leqslant
a<b\leqslant t_{k}$. If $S(t)$ has no dis\-con\-ti\-nu\-i\-ti\-es at
$(a,b)$, then the Riemann\,-von Mangoldt formula together with
Lagrange's mean value theorem imply that
\begin{multline}\label{lab_37}
S(b)\,-\,S(a)\,=\,(b-a)S'(c)\,=\,(b-a)\biggl(\,-\frac{1}{2\pi}\ln\frac{c}{2\pi}\,+\,o(1)\biggr)\,=\\
=\,-(b-a)\bigl(L_{k}\,+\,o(1)\bigr),\quad
L_{k}\,=\,\frac{1}{2\pi}\ln\frac{t_{k}}{2\pi}
\end{multline}
for some $c$, $a<c<b$. The relation (\ref{lab_37}) holds true if
$a$ or $b$ coincides with the ordinates of zeta zeros. In this cases, one should replace $S(a), S(b)$ by
$S(a+0), S(b-0)$, respectively.

Suppose that $\gamma_{(1)}<\ldots < \gamma_{(r)}$ are all the ordinates of zeros of $\zeta(s)$ lying on $[a,b]$,
and let $\kappa_{(1)},\ldots, \kappa_{(k)}$ be their multiplicities. Then we have:
\begin{multline}\label{lab_38}
S(b-0)\,-\,S(a+0)\,=\,\bigl(S(b-0)\,-\,S(\gamma_{(k)}+0)\bigr)\,+\,
\bigl(S(\gamma_{(k)}+0)-S(\gamma_{(k)}-0)\bigr)\,+\\
+\bigl(S(\gamma_{(k)}-0)-S(\gamma_{(k-1)}+0)\bigr)\,+\ldots\,+\,
\bigl(S(\gamma_{(1)}+0)-S(\gamma_{(1)}-0)\bigr)\,+\,\bigl(S(\gamma_{(1)}-0)-S(a+0)\bigr)\\
=\,\kappa_{(1)}\,+\ldots\,+\kappa_{(k)}\,-\,(b-a)\bigl(L_{k}+o(1)\bigr)\,\geqslant-\,(b-a)\bigl(L_{k}+o(1)\bigr)
\,\geqslant\\
\geqslant\,-\,(t_{k}-t_{k-1})\bigl(L_{k}+o(1)\bigr)\,=\,-1\,-\,o(1)
\end{multline}
(see Fig. 4).

Now we define $m$ and $n$ from the relations
$t_{m-1}<\tau-h\leqslant t_{m}$, $t_{n}\leqslant \tau+h<t_{n+1}$.
Suppose first that both points $\tau\pm h$ differs from the
ordinates of zeta zeros. By (\ref{lab_38}), we have:
\[
S(t_{m}-0)\,-\,S(\tau-h)\,\geqslant\,-1-o(1),\quad
S(\tau+h)\,-\,S(t_{n}+0)\,\geqslant\,-1-o(1),
\]
and hence
\begin{equation}\label{lab_39}
\Delta(m)\,=\,S(t_{m}+0)\,\geqslant\,S(t_{m}-0)\,\geqslant\,S(\tau-h)-1-o(1),
\end{equation}
\begin{equation}\label{lab_40}
\Delta(n)\,=\,S(t_{n}+0)\,\leqslant\,S(\tau+h)+1+o(1).
\end{equation}

Subtracting (\ref{lab_39}) from (\ref{lab_40}), we find:
\[
\Delta(n)\,-\,\Delta(m)\,\leqslant\,M+2+o(1)\,<\,M+3.
\]

Suppose now that $\tau + h$ is the ordinate of multiplicity
$\kappa\geqslant 1$. Then (\ref{lab_38}) implies
\[
S(\tau + h-0)\,-\,S(t_{n-1}+0)\,\geqslant\,-2-o(1),
\]
and therefore
\begin{equation}\label{lab_41}
\Delta(n-1)\,\leqslant\,S(\tau+h)+2+o(1)\,=\,S(\tau+h)\,-\,0.5\kappa\,+\,2+o(1)\,\leqslant\,S(\tau+h)\,+\,1.5\,+o(1).
\end{equation}

\begin{center}
\includegraphics{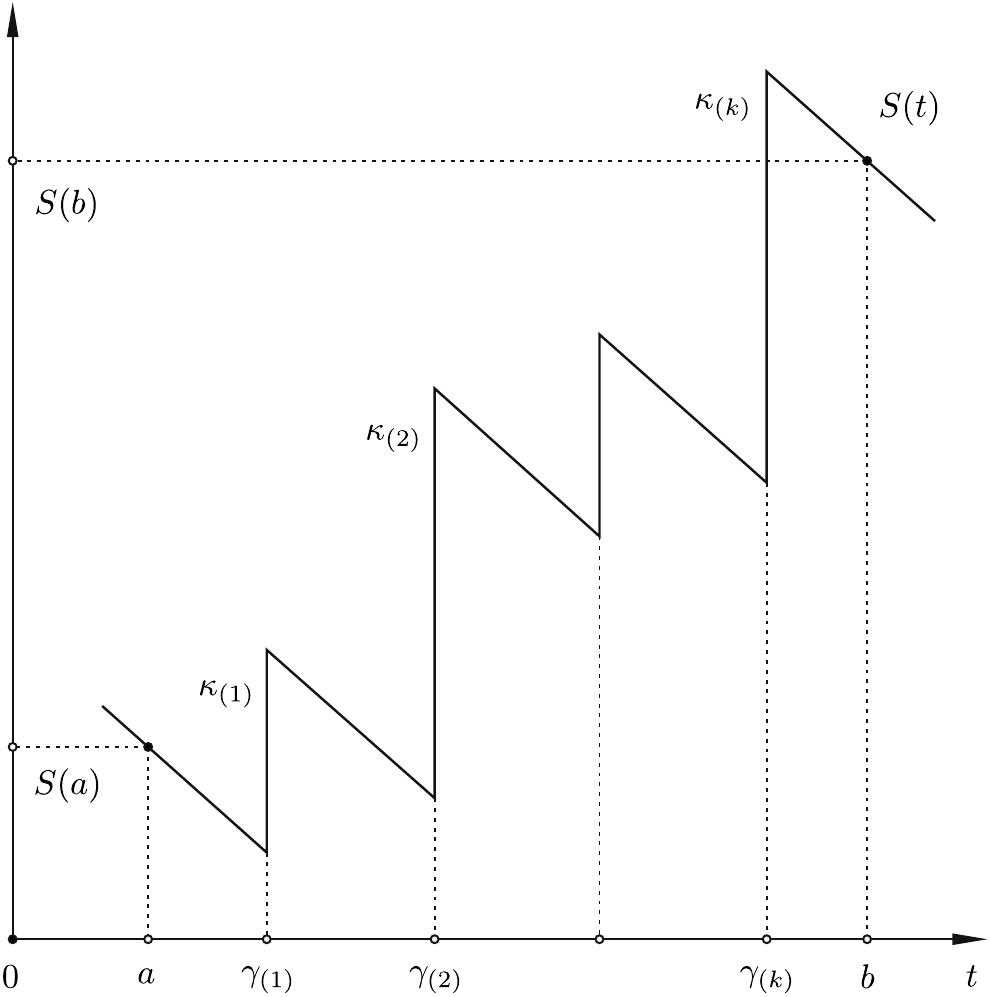}

\fontsize{10}{10pt}\selectfont Fig. 4. At each point $\gamma_{(r)}$ of discontinuity, the function $S(t)$ makes a jump equal
to the multiplicity of the ordinate $\gamma_{(r)}$, that is, to the sum of multiplicities of all zeta zeros with this point as ordinate.
\fontsize{12}{15pt}\selectfont
\end{center}

In view of (\ref{lab_39}), we get
\[
\Delta(n-1)\,-\,\Delta(m)\,\leqslant\,M+2.5+o(1)\,<\,M+3.
\]

Similarly, if $\tau - h$ is an ordinate of a zero of $\zeta(s)$, then
\[
S(t_{m+1}-0)\,-\,S(\tau - h)\,\geqslant\,-2-o(1),
\]
and hence
\begin{equation}\label{lab_42}
\Delta(m+1)\,=\,S(t_{m+1}+0)\,\geqslant\,
S(t_{m+1}-0)\,\geqslant\,S(\tau+h+0)-2-o(1)\,\geqslant\,S(\tau-h)-1.5-o(1).
\end{equation}
Taking (\ref{lab_40}) into account, we find
\[
\Delta(n)\,-\,\Delta(m+1)\,\leqslant\,M+2.5+o(1)\,<\,M+3.
\]
Finally, let both the points $\tau\pm h$ be the ordinates. By (\ref{lab_41}) and (\ref{lab_42}), we then have
\[
\Delta(n-1)\,-\,\Delta(m+1)\,\leqslant\,M+3+o(1)\,<\,M+3+10^{-4}.
\]
The above estimates imply that the smallest difference among
$\Delta(n-i)-\Delta(m+j)$, $0\leqslant i,j\leqslant 1$, does not
exceed $M+3+10^{-4}$ in any case. Denote by $n_{1}$ and $m_{1}$ the
corresponding values of $n-i$ and $m+j$ and set
$N=\bigl[-(M+3+10^{-4})\bigr]$. Since $N\geqslant 1$, we get
\begin{equation}\label{lab_43}
\bigl(\Delta(n_{1})\,-\,\Delta(n_{1}-1)\bigr)\,+\,\bigl(\Delta(n_{1}-1)\,-\,\Delta(n_{1}-2)\bigr)\,+\ldots
+ \bigl(\Delta(m_{1}+1)\,-\,\Delta(m_{1})\bigr)\,\leqslant\,-N.
\end{equation}
Formula (\ref{lab_36}) implies that $\Delta(k)-\Delta(k-1)\geqslant
-1$ and the equality takes place if and only if Gram's interval
$G_{k}$ is free of zeros of $\zeta\bigl(0.5+it\bigr)$. Thus,
(\ref{lab_43}) means that there are at least $N$ negative
differences (i.e. equal to $-1$) among $\Delta(k)-\Delta(k-1)$, $k =
m+1, \ldots, n$. Hence, there are at least $N$ intervals among the
intervals $G_{k}$, $k = m+1, \ldots, n$, which are free of zeros of
$\zeta\bigl(0.5+it\bigr)$.

To end the proof, we note that
\[
N\,\geqslant\,\frac{e^{a}}{10\sqrt{a}}\,-\,4\,>\,\frac{e^{a}}{16\sqrt{a}}\,=\,\frac{\sqrt{\pi
\varepsilon}}{16}\,\exp{\bigl((\pi\varepsilon)^{-1}\bigr)}\,
>\,0.1\sqrt{\varepsilon}\exp{\bigl((\pi\varepsilon)^{-1}\bigr)},
\]
and that all the intervals $G_{k}$, $k = m+1,\ldots, n$ are
contained in the segment $[\tau-h,\tau+h]$ of length $2h =
\varepsilon$. Theorem is proved.

\vspace{0.2cm}

The Corollary of Theorem 3 implies similar (but weaker) result for the distribution
of intervals $G_{n}$ containing at least two zeros of $\zeta(s)$.

\vspace{0.2cm}

\textsc{Theorem 6.} \emph{Suppose that the Riemann hypothesis is
true. Then there exist constants $T_{0}=T_{0}(\varepsilon)$ and
$c_{0} = c_{0}(\varepsilon)$ such that any segment $[T-H,T+H]$,
where $T\geqslant T_{0}$ and $H = 0.8\ln\ln\ln{T}+c_{0}$, contains
an interval $G_{k}$ with at least two zeros of $\zeta(s)$.}

\vspace{0.2cm}

\textsc{Proof.} By Corollary of Theorem 3, for sufficiently large $c$ and $H_{1} = 0.4\ln\ln\ln{T_{1}}+c$, the interval $(T_{1}-H_{1}, T_{1}+H_{1})$
contains a point $\tau_{1}$ such that $S(\tau_{1})<-3-10^{-3}$, and  the interval $(T_{1}+H_{1}, T_{1}+3H_{1})$ contains a point $\tau_{2}$ such that
$S(\tau_{2})>3+10^{-3}$.

We define $m$, $n$ by the inequalities $t_{m}<\tau_{1}\leqslant
t_{m+1}$, $t_{n-1}<\tau_{2}\leqslant t_{n}$. Using the same
arguments as in the proof of Theorem 4 together with the
inequalities $\tau_{1}<\tau_{2}$, $S(\tau_{2})-S(\tau_{1})>6+2\cdot
10^{-3}$, we find
\[
S(\tau_{1}-0)\,-\,S(t_{m}+0)\,\geqslant\,-1-o(1),
\]
and hence
\[
-\Delta(m)\,\geqslant\,-S(\tau_{1}-0)-1-o(1)\,\geqslant\,-S(\tau_{1})-1-o(1).
\]
Similarly,
\begin{multline*}
S(t_{n}+0)\,-\,S(\tau_{2})\,=\,\bigl(S(t_{n}+0)\,-\,S(t_{n}-0)\bigr)\,+\,\bigl(S(t_{n}-0)\,-\,S(\tau_{2}+0)\bigr)\,+\\
+\,\bigl(S(\tau_{2}+0)\,-\,S(\tau_{2})\bigr)\,\geqslant\,-1-o(1),
\end{multline*}
so we have $\Delta(n)\geqslant S(\tau_{2})-1-o(1)$. Therefore,
\[
\Delta(n)\,-\,\Delta(m)\,\geqslant\,S(\tau_{2})\,-\,S(\tau_{1})\,-\,2-o(1)\,>\,4.
\]
Thus, the inequality $\Delta(k)-\Delta(k-1)\geqslant 1$ holds for at
least one index $k$, $k = m+1,\ldots, n$. In view of (\ref{lab_36}),
the corresponding Gram's interval $G_{k}$ contains at least two
zeros of $\zeta\bigl(0.5+it\bigr)$. This interval lies in the
segment $[T_{1}-H_{1}, T_{1}+3H_{1}+t_{n}-t_{n-1}]$ whose length is
less than $1.6\ln\ln\ln{T_{1}}+4c+10^{-3}$. Setting $c_{0} =
2c+10^{-3}$, we arrive at the desired assertion. Theorem is proved.

\vspace{0.2cm}

Let $\gamma_{n}>0$ be an ordinate of a zero of $\zeta(s)$. Given
$n$, we indicate the unique number $m = m(n)$ such that
$t_{m-1}<\gamma_{n}\leqslant t_{m}$. Following Selberg
\cite{Selberg_1946b}, we denote $\Delta_{n} = m-n$. It is known (see
\cite[p. 355, remark 1]{Selberg_1989} and \cite{Korolev_2010}) that
$\Delta_{n}\ne 0$ for ``almost all'' $n$. Moreover, one can show
that the number of indices $n\leqslant N$ satisfying the condition
\[
\Delta_{n}\,\leqslant\,\frac{x}{\pi\sqrt{2}}\,\sqrt{\ln\ln{N}}
\]
is expressed as
\[
N\biggl(\;\frac{1}{\sqrt{2\pi}}\int_{-\infty}^{x}e^{-u^{2}/2}\,du\,+\,O\biggl(\frac{\ln\ln\ln{N}}{\sqrt{\ln\ln{N\mathstrut}}}\biggr)\biggr)
\]
for any $x$ (see \cite[Th. 5]{Boyarinov_2011c} and \cite[Th.
4-6]{Korolev_2012}). Given $N\geqslant N_{0}$, the above Theorem 3
allows to point out $M = M(N)$ such that the interval $(N,N+M]$
certainly contains an index $n$ with the condition $\Delta_{n}\ne
0$. Moreover, the following assertion holds.

\vspace{0.2cm}

\textsc{Theorem 7.} \emph{Suppose that the Riemann hypothesis is
true. Then there exist constants  $N_{0}$ and $c_{0} =
c_{0}(\varepsilon)$ such that the interval $(N,N+M]$, where
$N\geqslant N_{0}$ and
\[
M\,=\,\biggl[\frac{31}{5\pi}\,(\ln{N}+c_{0})\ln\ln\ln{N}\biggr],
\]
contains indices $n, m$ with the conditions $\Delta_{n} = 3$, $\Delta_{m}=-3$.}

\vspace{0.2cm}

\textsc{Proof.} We precede the proof by some remarks.

Firstly, the analogue of intermediate value theorem holds true for the function $S(t)$.
Namely, if $\tau_{1}<\tau_{2}$ and $S(\tau_{1})>S(\tau_{2})$ then for any $\alpha$ with the condition
$S(\tau_{2})<\alpha<S(\tau_{1})$, there exists a point $\tau$ on the interval $(\tau_{1}, \tau_{2})$ such that
$S(t)$ is continuous at this point and $S(\tau) =\alpha$ (see \cite[proof of Th. 3]{Korolev_2012b}).

Secondly, the value $S(t)$ is integer if and only if $t$ is Gram point (see \cite[proof of Th. 1]{Korolev_2012b}).

Suppose now that $T$ is sufficiently large. By Corollary of Theorem
3, for sufficiently large $c_{1}>0$ and $h = 0.4\ln\ln\ln{T}+c_{1}$,
the interval $(T,T+3h)$ contains the points $\tau_{1}<\tau_{2}$ such
that $S(\tau_{1})>3+10^{-3}$, $S(\tau_{2})<-3-10^{-3}$. By the first
remark, there exist a point $t$ between $\tau_{1}$ and $\tau_{2}$
such that $S(t)=S(t+0) = -3$. By second remark, this point is Gram
point, that is, $t = t_{\nu_{0}}$,
$S(t_{\nu_{0}}+0)=\Delta(\nu_{0})=-3$ for some $\nu_{0}$.

Similarly, we prove that each of intervals
$\bigl(T+(4j-1)h,T+(4j+3)h\bigr)$, $j = 1,\ldots,5$, contains Gram
point $t_{\nu_{j}}$ such that $S(t_{\nu_{j}}+0) =
\Delta(\nu_{j})=-3$. Now we take $T = t_{N}$. Since
\[
h = 0.4\ln\ln\ln{t_{N}}+c_{1}\,<\,0.4\ln\ln\ln{N}\,+\,c_{1},
\]
then the index $\nu$ defined by the relations
$t_{N+\nu}<T+23h\leqslant t_{N+\nu+1}$, satisfies the following
condition:
\[
\nu\,=\,
\frac{1}{\pi}\bigl(\vartheta(t_{N+\nu})\,-\,\vartheta(t_{N})\bigr)\,<\,\frac{23h}{\pi}\,\vartheta'(t_{N+\nu})\,<\,\frac{23h}{2\pi}\,\ln{N}\,<\,M.
\]
Hence, the interval $(N,N+\nu]$ contains at least $6$ indices $\nu_{j}$, $j = 0,\ldots, 5$, such that $\Delta(\nu_{j}) =-3$.
It is known (see \cite[Lemma 2]{Korolev_2012}) that the number of indices of the same interval satisfying the condition
$\Delta_{n} = 3$ differs from the above quantity for at most  $3 + (3-1) = 5$ in modulus. Hence, it is positive.

The proof of the second assertion of the Theorem is similar. It uses the fact that the difference between the number of indices
$n$ satisfying the condition $\Delta_{n} = -3$ and the number of indices with the condition $\Delta(\nu) = 3$
lying in the same interval, does not exceed $|-3|+|-3-1| = 7$. Theorem is proved.

\renewcommand{\refname}{\normalsize{Bibliography}}

\fontsize{11}{15pt}\selectfont

\noindent\textsc{Maxim A. Korolev}\\
Steklov Mathematical Institute\\
Russian Academy of Sciences\\
119991 Moscow, Russia\\
Gubkina st., 8;\\
E-mail: \texttt{korolevma@mi.ras.ru}\\
National Research Nuclear University ``MEPhI''\\
115409, Moscow, Russia\\
Kashirskoye shosse 31.
\end{document}